\input ./style/arxiv-general.cfg
\documentclass[aap,MSNbibl,nameyear,seceqn,secthm,dvips]{arximspdf}
\makeatletter
   \@ifpackageloaded{graphicx}{}{\usepackage{graphicx}}
\makeatother
\usepackage{mathrsfs}

\doi{10.1214/14-AAP1085}% Updated by VTEXPTS2LaTeX.exe, 05.02.2015
\volume{26}
\issue{1}
\pubyear{2016}
\firstpage{45}
\lastpage{72}
\docsubty{FLA}

\makeatletter
\newcommand{\perpp}{\perp\!\!\!\!\perp}
\renewcommand{\S}{\mathcal{S}}
\renewcommand{\i}{\mathrm{i}}
\newtheorem{lem}{Lemma}[section]
\newproclaim{example}{Example}[section]
\newtheorem{prop}{Proposition}[section]
\newproclaim{rem}{Remark}[section]
\newproclaim{assumption}{Assumption}
\makeatother

\begin{document}
\begin{frontmatter}

\title{A uniform law for convergence to the local times of linear fractional
stable motions}
\runtitle{Uniform convergence to local times}

\begin{aug}
% Corresponding author: James Duffy - james.duffy@economics.ox.ac.uk% Updated by VTEXPTS2LaTeX.exe, 08.02.2015 12:37
%Updated by VTEXPTS2LaTeX.exe, 06.02.2015 14:10
%Updated by VTEXPTS2LaTeX.exe, 05.02.2015 15:23
\author[A]{\fnms{James A.}~\snm{Duffy}\corref{}\ead[label=e1]{james.duffy@economics.ox.ac.uk}}%,
%\author[]{\fnms{}~\snm{}\ead[label=]{}}
% \and
%\author[]{\fnms{}~\snm{}\ead[label=]{}}
\runauthor{J. A. Duffy}
\affiliation{Institute for New Economic Thinking, Oxford Martin
School\\
and  University of Oxford}
%\dedicated{}
\address[A]{Nuffield College\\
Oxford\\
OX1 1NF\\
United Kingdom\\
\printead{e1}}
%\address[]{\\\printead{}}
\end{aug}

% HISTORY:
%
\received{\smonth{5} \syear{2014}}% Updated by VTEXPTS2LaTeX.exe,
%05.02.2015 15:23
%
\revised{\smonth{10} \syear{2014}}% Updated by VTEXPTS2LaTeX.exe,
%05.02.2015 15:23

% ABSTRACT
%
\begin{abstract}
We provide a uniform law for the weak convergence of additive functionals
of partial sum processes to the local times of linear fractional stable
motions, in a setting sufficiently general for statistical applications.
Our results are fundamental to the analysis of the global properties
of nonparametric estimators of nonlinear statistical models that involve
such processes as covariates.
\end{abstract}

% KEYWORDS
\end{frontmatter}

\section{Introduction}
\label{sec:intro}

Let $x_{t}=\sum_{s=1}^{t}v_{s}$ be the partial sum of a scalar linear
process $\{v_{t}\}$, for which the finite-dimensional distributions
of $d_{n}^{-1}x_{\lfloor nr\rfloor}$ converge to those of $X(r)$. Under
certain regularity conditions, we then have the finite-dimensional
convergence
\begin{equation}
\label{eq:fidiintro} \mathcal{L}_{n}^{f}(a,h_{n}):=
\frac{d_{n}}{nh_{n}}\sum_{t=1}^{n}f \biggl(
\frac
{x_{t}-d_{n}a}{h_{n}} \biggr)\mathop{\rightsquigarrow}_{\mathrm{f.d.d.}}\mathcal{L}(a)
\int_{\mathbb{R}}f,
\end{equation}
where $a\in\mathbb{R}$, $f$ is Lebesgue integrable, $h_{n}=o(d_{n})$
is a deterministic sequence, and $\mathcal{L}$ denotes the occupation
density (or \textit{local time}; see Remark~\ref{rem:locconv} below) associated
to $X$. Convergence results of this kind are particularly well documented
in the case where $\{x_{t}\}$ is a random walk [see the monograph
by \citet{BI95Stek}], and have more recently been extended to cover
generating mechanisms that allow the increments of $\{x_{t}\}$ to
exhibit significant temporal dependence [\citet{Jeg04,WP09ET}].

These more general theorems concerning (\ref{eq:fidiintro}) have, in
turn, played a fundamental role in the study of nonparametric estimation
and testing in the setting of nonlinear cointegrating models. The
simplest of these models takes the form
\begin{equation}
\label{eq:regress} y_{t}=m_{0}(x_{t})+u_{t},
\end{equation}
where $\{x_{t}\}$ is as above, $\{u_{t}\}$ is a weakly dependent
error process, and $m_{0}$ is an unknown function, assumed to possess
a certain degree of smoothness (or be otherwise approximable). In
a series of recent papers, (\ref{eq:fidiintro}) has facilitated the
development of a pointwise asymptotic distribution theory for kernel
regression estimators of $m_{0}$ under very general conditions: see
especially \citeauthor{WP09ET} (\citeyear{WP09ET,WP09Ecta,WP11ET,WP13mimeo}), \citet{KP12JoE}
and \citet{KAP12unpl}.\vspace*{2pt}%
\footnote{If $\{x_{t}\}$ is Markov, then this distribution theory may
be developed
by quite different arguments, without the use of (\ref{eq:fidiintro}),
see \citet{KT01AS} and \citet{KMT07AS}, which have spawned
a large
literature. While we consider this approach to the problem to be equally
important, our results touch upon it only a little, since we work
with a class of regressor processes that are typically (excepting
the random walk case) non-Markov.}

However, there are definite limits to the range of problems that can
be successfully addressed with the aid of (\ref{eq:fidiintro}). In particular,
since it concerns only the finite-dimensional convergence of $\mathcal{L}_{n}
^{f}(a,h_{n})$,
(\ref{eq:fidiintro}) is suited only to studying the \textit{local} behavior
of a nonparametric estimator: that is, its behavior in the vicinity
of a fixed spatial point. For the purpose of obtaining uniform rates
of convergence for kernel regression estimators on ``wide'' domains---that is, on domains having a width of the \textit{same} order as
the range of $\{x_{t}\}_{t=1}^{n}$---it is manifestly inadequate.
[See \citet{Duffy14rates}, for a detailed account.] The situation
is even worse with regard to sieve nonparametric estimation in this
setting---which initially motivated the author's research on this
problem---since in this case the development of even a \textit{pointwise}
asymptotic distribution theory requires a prior result on the \textit{uniform}
consistency of the estimator, over the entire domain on which estimation
is to be performed.

The main purpose of this paper is thus to provide conditions under
which the finite-dimensional convergence in (\ref{eq:fidiintro}) can
be strengthened to the weak convergence
\begin{equation}
\label{eq:unif} \mathcal{L}_{n}^{f}(a,h_{n})
\rightsquigarrow\mathcal{L}(a)\int_{\mathbb{R}}f,
\end{equation}
where $\mathcal{L}_{n}^{f}(a,h_{n})$ is regarded as a process indexed by
$(f,a)\in\mathscr{F}\times\mathbb{R}$, and $\{h_{n}\}$ may be
random. Results
of this kind are available in the existing literature, but only in
the random walk case, which requires that the increments of $\{x_{t}\}$
be independent, and $X$ to be an $\alpha$-stable L\'{e}vy motion [see
\citeauthor{borodin81TPA} (\citeyear{borodin81TPA,Bor82TP}); \citet{Perkins82PTRF}; and
\citet{BI95Stek}, Chapter~V].
In contrast, we allow the increments of $\{x_{t}\}$ to be serially
correlated, such that the associated limiting process $X$ may be
a linear fractional stable motion, which subsumes the $\alpha$-stable
L\'{e}vy motion and fractional Brownian motion as special cases. Further,
we permit the bandwidth sequence $\{h_{n}\}$ to be a \textit{random}
process, subject only to certain weak asymptotic growth conditions:
this is of considerable utility in statistical applications, where
the assumption that $\{h_{n}\}$ is a ``given'' deterministic sequence
seems quite unrealistic. Crucial to the proof of (\ref{eq:unif}) is
a novel order estimate for $\mathcal{L}_{n}^{f}(a,1)$ when $\int f=0$, which
is of interest in its own right.

The remainder of this paper is organized as follows. Our assumptions
on the data generating mechanism are described in Section~\ref{sec:assump}.
The main result (Theorem~\ref{thm:generalIP}) is discussed in
Section~\ref{sec:uniform}.
An outline of the proof follows in Section~\ref{sec:outline}, together with
the statement of two key auxiliary results (Propositions~\ref{prop:increment}
and \ref{prop:maxSn}). A preliminary application of our results to
the kernel nonparametric estimation of $m_{0}$ in (\ref{eq:regress})
is given in Section~\ref{sub:npreg}. The proof of Theorem~\ref{thm:generalIP} appears
in Section~\ref{sec:IPproof}, followed in Section~\ref{sec:ptwise} by
proofs of
Propositions~\ref{prop:increment}
and~\ref{prop:maxSn}. A proof related to the application appears
in Section~\ref{sec:proofsupp}. The final two Sections~\ref{sec:prelim}
and~\ref{sec:bndaux} are of a more technical nature, detailing the
proofs of two lemmas required in Section~\ref{sec:ptwise}, and so may
be skipped
on a first reading.

\subsection{Notation}
\label{sub:notation}

For a complete\vspace*{2pt} listing of the notation used in this paper, see Section~H
of the Supplement.%
\footnote{The Supplement is available as an addendum to \href{http://arxiv.org/abs/1501.05467}{arXiv:1501.05467}.}
The stochastic order notations $o_{p}(\cdot)$ and $O_{p}(\cdot)$
have the usual definitions, as given, for example, in \citet{VV98}, Section~2.2.
For deterministic sequences $\{a_{n}\}$ and $\{b_{n}\}$, we write
$a_{n}\sim b_{n}$ if $\lim_{n\rightarrow\infty}a_{n}/b_{n}=1$, and
$a_{n}\asymp b_{n}$
if $\lim_{n\rightarrow\infty}a_{n}/b_{n}\in(-\infty,\infty
)\setminus\{0\}$;
for random sequences, $a_{n}\lesssim_{p}b_{n}$ denotes $a_{n}=O_{p}(b_{n})$.
$X_{n}\rightsquigarrow X$ denotes weak convergence in the sense of
\citet{VVW96},
and $X_{n}\mathop{{\rightsquigarrow}}_{\mathrm{f.d.d.}} X$ the convergence of
finite-dimensional distributions.
For a metric space $(Q,d)$, $\ell_{\infty}(Q)$ [resp.,  $\ell
_{\mathrm{ucc}}(Q)$]
denotes the space of uniformly bounded functions on $Q$, equipped
with the topology of uniform convergence (resp., uniform convergence
on compacta). For $p\geq1$, $X$ a random variable, and $f\dvtx \mathbb{R}\rightarrow
\mathbb{R}$,
$\Vert X\Vert_{p}:=(\mathbb{E}\vert X\vert^{p})^{1/p}$ and $\Vert
f\Vert_{p}:=(\int_{\mathbb{R}}\vert f\vert^{p})^{1/p}$.
$\mathrm{BI}$ denotes the space of bounded and Lebesgue integrable functions
on $\mathbb{R}$. $\lfloor\cdot\rfloor$ and $\lceil\cdot\rceil$, respectively,
denote the floor and ceiling functions. $C$ denotes a generic constant
that may take different values even at different places in the same
proof; $a\lesssim b$ denotes $a\leq Cb$.

\section{Model and assumptions}
\label{sec:assump}

Our assumptions on the generating mechanism are similar to those of
\citet{Jeg04}---who proves a finite-dimen\-sional counterpart to our
main theorem---and are comparable to those made on the regressor
process in previous work on the estimation of nonlinear cointegrating
regressions [see, e.g.,  \citet{PP01Ecta}, \citeauthor{WP09Ecta} (\citeyear{WP09Ecta,WP12AS,WP13mimeo});
and \citet{KP12JoE}].

\begin{assumption}
\label{ass:reg}
\textup{(i)} $\{\varepsilon_{t}\}$ is a scalar i.i.d. sequence.
$\varepsilon_{0}$ lies in the domain of attraction of a strictly stable
distribution with index $\alpha\in(0,2]$, and has characteristic
function $\psi(\lambda):=\mathbb{E}\mathrm{e}^{\i\lambda\varepsilon
_{0}}$ satisfying
$\psi\in L^{p_{0}}$ for some $p_{0}\geq1$.

\textup{(ii)} $\{x_{t}\}$ is generated according to
\begin{equation}
\label{eq:regproc} x_{t} := \sum_{s=1}^{t}v_{s},
\qquad v_{t} :=\sum_{k=0}^{\infty
}
\phi_{k}\varepsilon_{t-k},
\end{equation}
and either:
\begin{longlist}[(a)]
\item[(a)] $\alpha\in(1,2]$, $\sum_{k=0}^{\infty}\vert
\phi _{k}\vert<\infty$
and $\phi:=\sum_{k=0}^{\infty}\phi_{k}\neq0$; or
$\phi_{k}\sim k^{H-1-1/\alpha}\pi_{k}$ for some $\{\pi_{k}\}_{k\geq0}$
strictly positive and slowly varying at infinity, with
\item[(b)] $H>1/\alpha$; or
\item[(c)] $H<1/\alpha$ and $\sum_{k=0}^{\infty}\phi_{k}=0$.
\end{longlist}

In both cases (b) and (c), $H\in(0,1)$.
\end{assumption}

\begin{rem}
\label{rem:levy}
Part (i) implies that there exists a
slowly varying sequence $\{\varrho_{k}\}$ such that
\begin{equation}
\label{eq:cvgtoZ} \frac{1}{n^{1/\alpha}\varrho_{n}}\sum_{t=1}^{\lfloor nr\rfloor
}
\varepsilon _{t}\mathop{\rightsquigarrow}_{\mathrm{f.d.d.}}Z_{\alpha}(r),
\end{equation}
where $Z_{\alpha}$ denotes an $\alpha$-\textit{stable L\'{e}vy motion}
on $\mathbb{R}$, with $Z_{\alpha}(0)=0$. That is,
the increments of $Z_{\alpha}$ are stationary, and for any $r_{1}<r_{2}$
the characteristic function of $Z_{\alpha}(r_{2})-Z_{\alpha}(r_{1})$
has the logarithm
\[
-(r_{2}-r_{1})c\vert\lambda\vert^{\alpha} \biggl[1+
\i\beta \operatorname{sgn}(\lambda)\tan \biggl(\frac{\pi\alpha}{2} \biggr) \biggr],
\]
where $\beta\in[-1,1]$ and $c>0$; following \citet{Jeg04}, page 1773,
we impose the further restriction that $\beta=0$ when $\alpha=1$.
We shall also require that $\{\varrho_{k}\}$ be chosen such that
$c=1$ here, which provides a convenient normalization for the scale
of $Z_{\alpha}$.
\end{rem}

\begin{rem}
To permit the alternative forms of (ii) to be more concisely
referenced, we shall \textit{refer} to (a) as corresponding
to the case where $H=1/\alpha$; this designation may be justified
by the manner in which the finite-dimensional limit of
$d_{n}^{-1}x_{\lfloor nr\rfloor}$
depends on $(H,\alpha)$, as displayed in (\ref{eq:Xdef}) below. The
statement that $H<1/\alpha$ will also be used as a shorthand for
(c), that is, it will always be understood that $\sum_{k=0}^{\infty}\phi_{k}=0$
in this case.
\end{rem}

We shall treat the parameters (including $H$ and $\alpha$) describing
the data generating mechanism as ``fixed'' throughout, ignoring the
dependence of any constants on these. Let $\{c_{k}\}$ denote a sequence
with $c_{0}=1$ and
\begin{equation}
\label{eq:ckdef} c_{k}= %
\cases{\phi, & \quad$\mbox{if }H=1/
\alpha$,\vspace*{3pt}
\cr
\vert H-1/\alpha\vert^{-1}k^{H-1/\alpha}
\pi_{k}, & \quad$\mbox {otherwise}$.}
\end{equation}
By Karamata's theorem [\citet{Bingham87}, Theorem~1.5.11], $\sum_{l=0}^{k}\phi_{k}\sim c_{k}$
as $k\rightarrow\infty$. Set
\begin{equation}
\label{eq:dkekdef} d_{k} := k^{1/\alpha}c_{k}
\varrho_{k}, \qquad e_{k} := kd_{k}^{-1},
\end{equation}
and note that the sequences $\{c_{k}\}$, $\{d_{k}\}$ and $\{e_{k}\}$
are regularly varying with indices $H-1/\alpha$, $H$ and $1-H$,
respectively. Theorems~5.1--5.3 in \citet{KM88Prob}
yield

\begin{prop}
\label{prop:fclt}
Under Assumption~\ref{ass:reg},
\begin{equation}
\label{eq:Xwkc} X_{n}(r):=\frac{1}{d_{n}}x_{\lfloor nr\rfloor}\mathop{
\rightsquigarrow }_{\mathrm{f.d.d.}}X(r),\qquad r\in [0,1],
\end{equation}
where $X$ is the linear fractional stable motion (LFSM)
\begin{eqnarray}
X(r) & :=& \int_{0}^{r}(r-s)^{H-1/\alpha}\,
\mathrm{d}Z_{\alpha
}(s)
\nonumber
\\[-8pt]
\label{eq:Xdef}
\\[-8pt]
\nonumber
&&{}+\int_{-\infty}^{0}
\bigl[(r-s)^{H-1/\alpha}-(-s)^{H-1/\alpha
}\bigr]\,\mathrm{d}Z_{\alpha}(s)
\end{eqnarray}
with the convention that $X=Z_{\alpha}$ when $H=1/\alpha$; $Z_{\alpha}$
is an $\alpha$-stable L\'{e}vy motion on $\mathbb{R}$, with
$Z_{\alpha}(0)=0$.
\end{prop}

\begin{rem}
For a detailed discussion of the LFSM, see \citet{ST94}. When
$\alpha=2$,
$Z_{\alpha}$ is a Brownian motion with variance $2$; if additionally
$H\neq1/\alpha$, $X$ is thus a fractional Brownian motion.
\end{rem}

\begin{rem}
\label{rem:cvgtypes}
Excepting such cases as the following:
\begin{longlist}[(iii)]
\item[(i)] $\alpha\in(1,2]$, $H>1/\alpha$ [\citet{Astr83Lith}, Theorem~2];
\item[(ii)] $\alpha=2$, $H=1/\alpha$ and $\mathbb{E}\varepsilon
_{0}^{2}<\infty$
[\citet{Hannan1979SP}];
and
\item[(iii)] $\alpha=2$, $H<1/\alpha$ and $\mathbb{E}\vert\varepsilon
_{0}\vert^{q}<\infty$
for some $q>2$ [\citet{DdJ00ET}, Theorem~3.1];
\end{longlist}
it may not be possible to strengthen the convergence in (\ref{eq:Xwkc})
to weak convergence on $\ell_{\infty}[0,1]$. Weak convergence may
hold, however, with respect to a weaker topology, and we shall be
principally concerned with whether this topology is sufficiently strong
that
\begin{equation}
\label{eq:extremecvg} \inf_{r\in[0,1]}X_{n}(r) \rightsquigarrow
\inf_{r\in[0,1]}X(r),\qquad \sup_{r\in
[0,1]}X_{n}(r)
\rightsquigarrow\sup_{r\in[0,1]}X(r),
\end{equation}
such as would follow from weak convergence in the Skorokhod $M_{1}$
topology [see \citet{Skor56TP}, Section~2.2.10]. When $H=1/\alpha$, sufficient
conditions for this kind of convergence---which entail further restrictions
on $\{\phi_{k}\}$ than are imposed here---are given in \citet{AT92AP}, Theorem~2
and \citet{Tyran10SPLett}, Theorem~1 and Corollary~1.
However, when
$H<1/\alpha$
and $\alpha\in(0,2)$, the sample paths of $X$ are unbounded, and
thus (\ref{eq:extremecvg}) cannot possibly hold [see \citet{ST94},
Example~10.2.5]. In any case, (\ref{eq:extremecvg}) is \textit{not} necessary
for the main results of this paper; it merely permits Theorem~\ref
{thm:generalIP}
below to take a slightly strengthened form.
\end{rem}

\begin{rem}
\label{rem:locconv}
In consequence of Theorem~3(i) in \citet{Jeg04},
the convergence in (\ref{eq:Xwkc}) occurs jointly with
\[
\mathcal{L}_{n}^{f}(a):=\frac{1}{e_{n}}\sum
_{t=1}^{n}f(x_{t}-d_{n}a)\mathop{\rightsquigarrow}_{\mathrm{f.d.d.}} \mathcal{L}(a)\int_{\mathbb{R}}f,
\qquad a\in\mathbb{R}
\]
for every $f\in\mathrm{BI}$. Here, $\{\mathcal{L}(a)\}_{a\in\mathbb
{R}}$ denotes
the occupation density (\textit{local time}) of~$X$, a process which,
almost surely, has continuous paths and satisfies
\begin{equation}
\label{eq:otf} \int_{\mathbb{R}}f(x)\mathcal{L}(x)\,\mathrm{d}x=\int
_{0}^{1}f\bigl(X(r)\bigr)\,\mathrm{d}r
\end{equation}
for all Borel measurable and locally integrable $f$. (For the existence
of $\mathcal{L}$, see Theorem~0 in \citet{Jeg04}; the path continuity
may be deduced from Theorem~\ref{thm:generalIP} below.)
\end{rem}

\section{A uniform law for the convergence to local time}
\label{sec:uniform}

Our main result concerns the convergence
\begin{equation}
\label{eq:mainIP} \mathcal{L}_{n}^{f}(a,h_{n}):=
\frac{1}{e_{n}h_{n}}\sum_{t=1}^{n}f \biggl(
\frac{x_{t}-d_{n}a}{h_{n}} \biggr)\rightsquigarrow\mathcal {L}(a)\int_{\mathbb{R}}f,
\end{equation}
where $\mathcal{L}_{n}^{f}(a,h_{n})$ is regarded as a process indexed by
$(f,a)\in\mathscr{F}\times\mathbb{R}$. ($\mathscr{F}\times\mathbb
{R}$ is endowed with
the product topology, $\mathscr{F}$ having the $L^{1}$ topology, and
$\mathbb{R}$
the usual Euclidean topology.) $\{h_{n}\}$ is a measurable bandwidth
sequence that may be functionally dependent on $\{x_{t}\}$, or indeed
upon any other elements of the probability space; it is required only
to satisfy:

\begin{assumption}
\label{ass:band}
$h_{n}\in\mathscr{H}_{n}:=[\underline
{h}_{n},\overline{h}_{n}]$
with probability approaching 1\break (w.p.a.1), where $\overline{h}_{n}=o(d_{n})$
and $\underline{h}_{n}^{-1}=o(e_{n}\log^{-2}n)$.
\end{assumption}

Define
\begin{equation}
\label{eq:BImom} \mathrm{BI}_{\beta}:= \biggl\{ f\in\mathrm{BI}\vert\int
_{\mathbb
{R}}\bigl\vert f(x)\bigr\vert\vert x\vert^{\beta}\,\mathrm{d}x<
\infty \biggr\}
\end{equation}
and let $\mathrm{BIL}_{\beta}$ denote the subset of Lipschitz continuous
functions in $\mathrm{BI}_{\beta}$. In order to state conditions on
$\mathscr{F}
\subset\mathrm{BI}$
that are sufficient for (\ref{eq:mainIP}) to hold, we first recall some
definitions familiar from the theory of empirical processes. A function
$F\dvtx \mathbb{R}\rightarrow\mathbb{R}_{+}$ is termed an \textit{envelope} for $\mathscr{F}$,
if $\sup_{f\in\mathscr{F}}\vert f(x)\vert\leq F(x)$ for every $x\in
\mathbb{R}$.
Given a pair of functions $l,u\in L^{1}$, define the bracket
\[
[l,u]:=\bigl\{f\in L^{1}\vert l(x)\leq f(x)\leq u(x), \forall x\in
\mathbb{R} \bigr\};
\]
we say that $[l,u]$ is an $\varepsilon$-\textit{bracket} if $\Vert
u-l\Vert_{1}<\varepsilon$,
and a \textit{continuous} bracket if $l$ and $u$ are continuous. Let
$N_{[\,]}^{\ast}(\varepsilon,\mathscr{F})$ denote the minimum number of\vspace*{2pt}
continuous
$\varepsilon$-brackets required to cover $\mathscr{F}$.

\begin{assumption}
\label{ass:fset}
\textup{(i)} $\mathscr{F}\subset\mathrm{BI}$ has
envelope $F\in\mathrm{BIL}_{\beta}$,
for some $\beta>0$; and

\textup{(ii)} for each $\varepsilon>0$, $N_{[\,]}^{\ast
}(\varepsilon,\mathscr{F})<\infty$.
\end{assumption}

We may now state our main result, the proof of which appears in
Section~\ref{sec:IPproof}.
\begin{thm}
\label{thm:generalIP}
Suppose Assumptions~\ref{ass:reg}--\ref{ass:fset}
hold. Then:
\begin{longlist}[(ii)]
\item[(i)] (\ref{eq:mainIP}) holds in $\ell_{\mathrm{ucc}}(\mathscr{F}\times
\mathbb{R})$;

\noindent and if additionally (\ref{eq:extremecvg}) holds, then

\item[(ii)] (\ref{eq:mainIP}) holds in $\ell_{\infty
}(\mathscr{F}\times\mathbb{R})$.
\end{longlist}
\end{thm}

\begin{rem}
The case where $h_{n}=1$, $\mathscr{F}=\{f\}$ and $\{x_{t}\}$ is a random
walk---which here entails $H=1/\alpha$ and $\phi_{i}=0$ for all
$i\geq1$---has been studied extensively: see in particular
\citeauthor{borodin81TPA} (\citeyear{borodin81TPA,Bor82TP}),
\citet{Perkins82PTRF} and \citet{BI95Stek}, Chapter~V.
In those works,
it is proved (under these more restrictive assumptions on $\{x_{t}\}$)
that
\[
\frac{1}{e_{n}}\sum_{t=1}^{\lfloor nr\rfloor
}f(x_{t}-d_{n}a)
\rightsquigarrow\mathcal{L} (a;r)\int_{\mathbb{R}}f
\]
on $\ell_{\infty}(\mathbb{R}\times[0,1])$, where $\mathcal
{L}(a;r)$ denotes
the local time of $X$ restricted to $[0,r]$. Theorem~\ref{thm:generalIP}
could be very easily extended in this direction; we have refrained
from doing so only to keep the paper to a reasonable length. The principal
contribution of Theorem~\ref{thm:generalIP} is thus to extend this convergence
in a direction more suitable for statistical applications, by allowing
$\{v_{t}\}$ to be serially correlated and the bandwidth sequence
$\{h_{n}\}$ to be data-dependent.
\end{rem}

\begin{rem}
\label{rem:LCN}
After the manuscript of this paper had been completed,
we obtained a copy of an unpublished manuscript by \citet{LCW14mimeo}
in which, under rather different assumptions from those given here,
a result similar to Theorem~\ref{thm:generalIP} is proved (for a fixed $f$
and a deterministic sequence $\{h_{n}\}$). Regarding the differences
between our main result and their Theorem~2.1, we may note particularly
their requirement that there exist a sequence of processes $\{
X_{n}^{\ast}\}$
with $X_{n}^{\ast}=_{d}X$, and a $\delta>0$ such that
\begin{equation}
\label{eq:LCWcondition} \sup_{r\in[0,1]}\bigl\vert X_{n}(r)-X_{n}^{\ast}(r)
\bigr\vert=o_{\mathrm
{a.s.}}\bigl(n^{-\delta
}\bigr),
\end{equation}
a condition which excludes a large portion of the processes considered
in this paper, in view of Remark~\ref{rem:cvgtypes} above. [The availability
of (\ref{eq:LCWcondition}) permits these authors to prove their result
by an argument radically different from that developed here.] On the
other hand, our results do \textit{not} subsume theirs, since these
authors do not require $v_{t}$ to be a linear process.
\end{rem}

Although Assumption~\ref{ass:fset} requires that $\mathscr{F}$ have a
smooth envelope
and smooth brackets, it is perfectly consistent with $\mathscr{F}$ containing
discontinuous functions. Indeed, Assumption~\ref{ass:fset} is
consistent with
such cases as the following, as verified in Section~A of the Supplement.
[We expect that boundedness and $\int\vert f(x)\vert\vert x\vert
^{\beta
}\,\mathrm{d}x<\infty$
could also be relaxed through the use of a suitable truncation argument,
such as is employed in the proof of Theorem~V.4.1 in
\citet{BI95Stek}.]

\begin{example}[(Single function)]
\label{exa:single}$\mathscr{F}=\{f\}$ where $f\in\mathrm{BI}_{\beta
}$, and
is majorised by another function $F\in\mathrm{BIL}_{\beta}$, in the sense
that $\vert f(x)\vert\leq F(x)$ for all $x\in\mathbb{R}$. This obtains
trivially if $f$ is itself in $\mathrm{BIL}_{\beta}$ (simply take $F(x):=
\vert f(x)\vert$),
but is also consistent with $f\in\mathrm{BI}_{\beta}$ having
finitely many
discontinuities (at the points $\{a_{k}\}_{k=1}^{K}$, where $a_{k}<a_{k+1}$),
and being Lipschitz continuous on $(-\infty,a_{1})\cup[a_{K},\infty)$;
all that is really necessary here is for $f$ to have one-sided Lipschitz
approximants. Importantly, this includes the case where $f(x)=\mathbf
{1}\{
x\in I\}$
for any bounded interval~$I$.
\end{example}

\begin{example}[(Parametric family)]
\label{exa:param}$\mathscr{F}=\{g(x,\theta)\vert\theta\in\Theta\}
\subset\mathrm{BIL}_{\beta}$,
where $\Theta$ is compact, and there exists a $\tau\in(0,1]$ and
a $\dot{g}\in\mathrm{BIL}_{\beta}$ such that
\[
\bigl\vert g(x,\theta)-g\bigl(x,\theta^{\prime}\bigr)\bigr\vert\leq\dot{g}(x)\bigl\Vert
\theta -\theta^{\prime}\bigr\Vert^{\tau}
\]
for all $\theta,\theta^{\prime}\in\Theta$.
\end{example}

\begin{example}[(Smooth functions)]
\label{exa:smooth} $\mathscr{F}=\{f\in C^{\tau}(\mathbb{R})\vert
\vert f\vert\leq F\}$,
where $F\in\mathrm{BIL}_{\beta}$ and
\[
C_{L}^{\tau}(\mathbb{R}):=\bigl\{f\in\mathrm{BI}\vert \ \exists
C_{f}<L\mbox{ s.t. }\bigl\vert f(x)-f\bigl(x^{\prime}\bigr)\bigr\vert\leq
C_{f}\bigl\vert x-x^{\prime}\bigr\vert ^{\tau} \ \forall
x,x^{\prime}\in\mathbb{R}\bigr\}
\]
for some $\tau\in(0,1]$ and $L<\infty$.
\end{example}

\section{Outline of proof and auxiliary results}
\label{sec:outline}

\subsection{Outline of proof}

\begin{figure}

\includegraphics{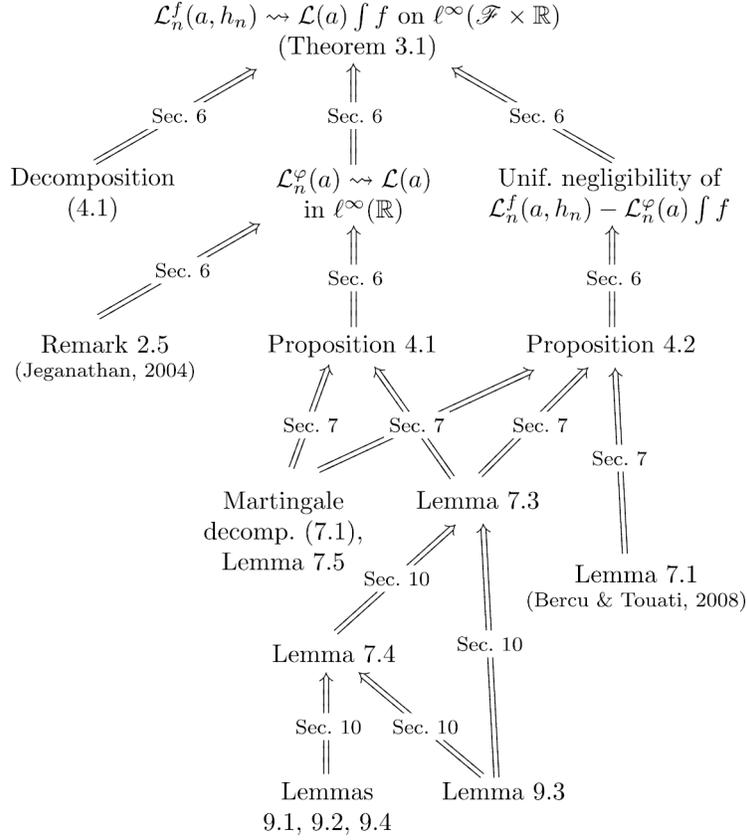}

\caption{Outline of proofs.}
\label{fig:mainrel}
\end{figure}

The principal relationships between the results in this paper are
summarized in Figure~\ref{fig:mainrel}. The proof of Theorem~\ref
{thm:generalIP},
depicted\vadjust{\goodbreak} in the top half of the figure, proceeds as follows. To reduce
the difficulties arising by the randomness of $h=h_{n}$, we decompose
\begin{equation}
\label{eq:Lndecomp} \mathcal{L}_{n}^{f}(a,h)=\mathcal{L}_{n}^{\varphi}(a)
\int_{\mathbb
{R}}f+ \biggl[\mathcal{L}_{n}
^{f}(a,h)-\mathcal{L}_{n}^{\varphi}(a)\int
_{\mathbb{R}}f \biggr],
\end{equation}
where
\begin{equation}
\label{eq:cmpct}
\varphi(x):=\bigl(1-\vert x\vert\bigr)\mathbf{1}\bigl\{\vert x\vert\leq1\bigr\}
\end{equation}
denotes the triangular kernel function, and $\mathcal{L}_{n}^{\varphi}(a):=
\mathcal{L}_{n}^{\varphi}(a,1)$.
(This choice of $\varphi$ is made purely for convenience; any compactly
supported Lipschitz function would serve our purposes equally well
here.) It thus suffices to show that $\mathcal{L}_{n}^{\varphi
}\rightsquigarrow\mathcal{L}$
in $\ell^{\infty}(\mathbb{R})$, and that the bracketed term on the right-hand
side of (\ref{eq:Lndecomp}) is uniformly negligible over $(f,h)\in
\mathscr{F}
\times\mathscr{H}_{n}$.

%f1
%

In view of Remark~\ref{rem:locconv} above, the finite-dimensional distributions
of $\mathcal{L}_{n}^{\varphi}$ converge to those of $\mathcal{L}$.
The asymptotic
tightness of $\mathcal{L}_{n}^{\varphi}$ will follow from the bound
on the
spatial increments
\[
\mathcal{L}_{n}^{\varphi}(a_{1})-\mathcal{L}_{n}^{\varphi
}(a_{2})=
\frac{1}{e_{n}}\sum_{t=1}^{n}\bigl[
\varphi(x_{t}-d_{n}a_{1})-\varphi
(x_{t}-d_{n}a_{2})\bigr]=:\frac
{1}{e_{n}}\sum
_{t=1}^{n}g_{1}(x_{t}),
\]
given in Proposition~\ref{prop:increment} below. The bracketed term on
the right-hand
side of (\ref{eq:Lndecomp}) may be written as
\begin{equation}
\label{eq:centring} \frac{1}{e_{n}}\sum_{t=1}^{n}
\biggl[\frac{1}{h}f \biggl(\frac
{x_{t}-d_{n}a}{h} \biggr)-
\varphi(x_{t}-d_{n}a)\int_{\mathbb
{R}}f \biggr]=:
\frac{1}{e_{n}}\sum_{t=1}^{n}g_{2}(x_{t}).
\end{equation}
Control of (\ref{eq:centring}) over progressively denser subsets of
$\mathscr{F}\times\mathscr{H}_{n}$ is provided by Proposition~\ref
{prop:maxSn} below; the
conjunction of a bracketing argument and the continuity of the brackets
suffices to extend this to the entirety of $\mathscr{F}\times\mathscr
{H}_{n}$.

By\vspace*{1pt} construction, both $\int g_{1}=0$ and $\int g_{2}=0$. The proofs
of Propositions\break \ref{prop:increment}~and~\ref{prop:maxSn} may therefore
be approached in a unified way, through the analysis of sums of the
form
\begin{equation}
\label{eq:Sng} \S_{n}g:=\sum_{t=1}^{n}g(x_{t}),
\end{equation}
where $g$ ranges over a class $\mathscr{G}$, all members of which have
the property that\break $\int g=0$. Such functions are termed \textit{zero
energy} functions [\citet{WP11ET}]; we shall correspondingly\vspace*{1pt} term
$\{\S_{n}g\}_{g\in\mathscr{G}}$ a \textit{zero energy process}. Such processes
are ``centered'' in the sense that $e_{n}^{-1/2}\S_{n}g$ converges
weakly to a mixed Gaussian variate [\citet{Jeg08CDFP}, Theorem~5];
whereas $e_{n}^{-1}\S_{n}g\rightsquigarrow\mathcal{L}(0)\int g$ if
$\int g\neq0$.

Equation (\ref{eq:Sng}) will be handled by decomposing $\S_{n}g$ as
\[
\S_{n}g=\sum_{k=0}^{n-1}
\mathcal{M}_{nk}g+\mathcal{N}_{n}g,
\]
where each $\mathcal{M}_{nk}g$ is a martingale; see (\ref{eq:SnDecomp}) below.
We provide order estimates for the sums of squares and conditional
variances of the $\mathcal{M}_{nk}g$'s (Lemma~\ref{lem:fundamental});
by an application
of either Burkholder's inequality, or a tail bound due to \citet{BT08AAP},
these translate into estimates for the $\mathcal{M}_{nk}g$'s themselves.
Propositions~\ref{prop:increment} and \ref{prop:maxSn} then follow
by standard arguments.

\subsection{Key auxiliary results}
\label{sub:zeroen}

To state these, we introduce the quantity
\begin{equation}
\label{eq:alphnorm} \Vert f\Vert_{[\beta]}:=\inf\bigl\{c\in\mathbb{R}_{+}
\vert\bigl\vert\hat {f}(\lambda )\bigr\vert\leq c\vert\lambda\vert^{\beta}, \forall
\lambda\in\mathbb{R}\bigr\}
\end{equation}
for $f\in\mathrm{BI}$, $\beta\in(0,1]$, and $\hat{f}(\lambda
):=\int\mathrm{e}^{\i
\lambda x}f(x)\,\mathrm{d}x$.
It is easily verified that $\Vert f\Vert_{[\beta]}$ is indeed a norm on
the space $\mathrm{BI}_{[\beta]}:=\{f\in\mathrm{BI}\vert\Vert
f\Vert_{[\beta ]}<\infty\}$
(modulo equality almost everywhere). Some useful properties of
$\Vert f\Vert_{[\beta]}$
are collected in Lemma~\ref{lem:alphnorm} below; in particular, it is shown
that $\mathrm{BI}_{[\beta]}$ contains all $f\in\mathrm{BI}_{\beta
}$ for which
$\int f=0$. Define
\begin{equation}
\label{eq:abvbeta} \overline{\beta}_{H}:=\frac{1-H}{2H}\wedge1,
\end{equation}
noting that $\overline{\beta}_{H}\in(0,1]$ for all $H\in(0,1)$, and let
$\Vert\cdot\Vert_{\tau_{2/3}}$ denote the Orlicz norm associated
to the convex and increasing function
\begin{equation}
\label{eq:tau23} \tau_{2/3}(x):= %
\cases{ x(e-1), & \quad$
\mbox{if }x\in[0,1]$,\vspace*{3pt}
\cr
\mathrm{e}^{x^{2/3}}-1, & \quad$
\mbox{if }x\in(1,\infty)$.}
\end{equation}
[See \citet{VVW96}, page~95 for the definition of an Orlicz norm.]
A bound on the spatial increments of $\mathcal{L}_{n}^{\varphi}$ is given
by:

\begin{prop}
\label{prop:increment}
For every $\beta\in(0,\overline{\beta}_{H})$,
there exists $C_{\beta}<\infty$ such that
\[
\sup_{a_{1},a_{2}\in\mathbb{R}}\bigl\Vert\mathcal{L}_{n}^{\varphi
}(a_{1})-
\mathcal{L}_{n}^{\varphi}(a_{2})\bigr\Vert_{\tau
_{2/3}}
\leq C_{\beta}\vert a_{1}-a_{2}\vert^{\beta}.
\]
\end{prop}

The next result shall be applied to prove that the recentered sums
(\ref{eq:centring}) are uniformly negligible. Since the order estimate
given below is of interest in its own right [see \citet{Duffy14rates},
for an example of how it may be used to determine the uniform order
of the first-order bias of a nonparametric regression estimator], we
shall state it at a slightly higher level of generality than is
needed for our purposes here. For $\mathscr{F}\subset\mathrm
{BI}_{[\beta]}$, define
\begin{equation}
\label{eq:delnbetaf} \delta_{n}(\beta,\mathscr{F}):=\Vert\mathscr{F}
\Vert_{\infty}+e _{n}^{1/2}\bigl(\Vert\mathscr{F}
\Vert_{1}+\Vert\mathscr{F}\Vert_{2}\bigr)+e _{n}d_{n}^{-\beta}
\Vert\mathscr{F}\Vert_{[\beta]},
\end{equation}
where $\Vert\mathscr{F}\Vert:=\sup_{f\in\mathscr{F}}\Vert f\Vert$.

\begin{prop}
\label{prop:maxSn}
Suppose $\beta\in(0,\overline{\beta}_{H})$ and
$\mathscr{F}
_{n}\subset\mathrm{BI}_{[\beta]}$
with $\#\mathscr{F}_{n}\lesssim n^{C}$. Then
\begin{equation}
\label{eq:Snbnd} \max_{f\in\mathscr{F}_{n}}\vert\S_{n}f\vert
\lesssim_{p}\delta _{n}(\beta,\mathscr{F} _{n})\log
n.
\end{equation}
If also $\Vert\mathscr{F}_{n}\Vert_{1}\lesssim1$, $\Vert\mathscr
{F}_{n}\Vert_{[\beta ]}=o(d_{n}^{\beta})$
and $\Vert\mathscr{F}_{n}\Vert_{\infty}=o(e_{n}\log^{-2}n)$, then
\[
\max_{f\in\mathscr{F}_{n}}\vert\S_{n}f\vert=o_{p}(e_{n}).
\]
\end{prop}

\begin{rem}
\label{rem:maxSnb1}
As is clear from the proof, if $\beta\in
[\overline{\beta}_{H},1]$
then (\ref{eq:Snbnd}) holds in a modified form, with $e
_{n}d_{n}^{-\beta}\Vert\mathscr{F}\Vert_{[\beta]}$
in (\ref{eq:delnbetaf}) being replaced by
\[
\Biggl[\sum_{k=1}^{n-1}d_{k}^{-(1+\beta)}+e_{n}^{1/2}
\sum_{k=1}^{n-1}k^{-1/2}d_{k}^{-(1+2\beta)/2}
\Biggr]\Vert\mathscr {F}\Vert_{[\beta]}.
\]
\end{rem}

The proofs of Propositions~\ref{prop:increment} and \ref{prop:maxSn}
are given in Section~\ref{sec:ptwise} below.

\section{A preliminary application to nonparametric regression}\label{sub:applic}\label{sub:npreg}

Suppose that we observe $\{(y_{t},x_{t})\}_{t=1}^{n}$ generated according
to the nonlinear cointegrating regression model
\[
y_{t}=m_{0}(x_{t})+u_{t},
\]
where $\{u_{t}\}$ is some weakly dependent disturbance process. As
shown in \citet{WP09Ecta}, under suitable smoothness conditions the
unknown function~$m_{0}$ may be consistently estimated, at each fixed
$x\in\mathbb{R}$, by the Nadaraya--Watson estimator
\begin{eqnarray}
\hat{m}(x) & =&\frac{\sum_{t=1}^{n}K_{h_{n}}(x_{t}-x)y_{t}}{\sum_{t=1}^{n}K_{h_{n}}(x_{t}-x)}
\nonumber
\\[-8pt]
\label{eq:kernreg}
\\[-8pt]
\nonumber
& =& m_{0}(x)+\frac{\sum_{t=1}^{n}K_{h_{n}}(x_{t}-x)[(m_{0}(x_{t})-m_{0}(x))+u_{t}]}{\sum_{t=1}^{n}K_{h_{n}}(x_{t}-x)},
\end{eqnarray}
where $K_{h}(u):=h^{-1}K(h^{-1}u)$, and $K\in\mathrm{BI}$ is a positive,
mean-zero kernel with $\int_{\mathbb{R}}K=1$.

Now consider the problem of determining the rate at which $\hat{m}$
converges uniformly to $m_{0}$. As a first step, we would need to
obtain the uniform rate of divergence of the denominator in (\ref{eq:kernreg});
it is precisely this rate that the preceding results allow us to compute.
By Theorem~\ref{thm:generalIP},
\begin{equation}
\label{eq:kernloccvg} \frac{1}{e_{n}}\sum_{t=1}^{n}K_{h_{n}}(x_{t}-d_{n}a)
\rightsquigarrow \mathcal{L}(a)
\end{equation}
in $\ell_{\infty}(\mathbb{R})$, provided $\mathbb{P}\{h_{n}\in
\mathscr{H}_{n}\}\rightarrow1$,
and $K$ satisfies the requirements of Example~\ref{exa:single} (which seems
broad enough to cover any reasonable choice of $K$). Since $a\mapsto
\mathcal{L}(a)$
is random---being dependent on the trajectory of the limiting process
$X$---we now face the problem of identifying a sequence of sets
on which the left side of (\ref{eq:kernloccvg}) can be uniformly bounded
away from zero. A natural candidate is
\[
A_{n}^{\varepsilon}:=\bigl\{x\in\mathbb{R}\vert\mathcal
{L}_{n}\bigl(d_{n}^{-1}x\bigr)\geq\varepsilon\bigr
\},
\]
where $\varepsilon>0$. $\mathcal{L}_{n}$ is trivially bounded away from zero
on this set, whence
\begin{equation}
\label{eq:denomrate} \sup_{x\in A_{n}^{\varepsilon}} \Biggl[\sum
_{t=1}^{n}K_{h_{n}}(x_{t}-x)
\Biggr]^{-1}\lesssim_{p}e_{n}^{-1}.
\end{equation}
More significantly, for any given $\delta>0$, we may choose $\varepsilon>0$
such that
\begin{equation}
\label{eq:mostsupp} \limsup_{n\rightarrow\infty}\mathbb{P} \Biggl\{
\frac{1}{n}\sum_{t=1}^{n}\mathbf{1}
\bigl\{ x_{t}\notin A_{n}^{\varepsilon}\bigr\}\geq\delta
\Biggr\} \leq\delta;
\end{equation}
see Section~\ref{sec:proofsupp} for the proof. That is, $\varepsilon>0$ may
be chosen such that $A_{n}^{\varepsilon}$ contains as large a fraction
of the observed trajectory $\{x_{t}\}_{t=1}^{n}$ as is desired, in
the limit as $n\rightarrow\infty$. [Were we to allow $\varepsilon
=\varepsilon
_{n}\rightarrow0$,
we could permit $\delta=\delta_{n}\rightarrow0$ here, but the order of
(\ref{eq:denomrate}) would necessarily be increased.]

Note that the sample-dependence of $A_{n}^{\varepsilon}$ is necessary
for it fulfill two roles here, by being both ``small'' enough for (\ref{eq:denomrate})
to hold, but also ``large'' enough to be consistent with (\ref{eq:mostsupp}).
If $A_{n}^{\varepsilon}$ were replaced by a sequence of deterministic
intervals (or sets, more generally), then the maintenance of (\ref{eq:denomrate})
would necessarily come at the cost of violating (\ref{eq:mostsupp}).
For example, when $x_{t}$ is a random walk with finite variance
($\alpha=2$),
the ``widest'' sequence of intervals $[-a_{n},a_{n}]$ for which (\ref{eq:denomrate})
holds is one for which $a_{n}=o(n^{1/2})$: but in consequence, the
\textit{fraction} of any trajectory $\{x_{t}\}_{t=1}^{n}$ falling within
such an interval will converge to $0$, as $n\rightarrow\infty$ [see
Remark~2.8 in \citet{Duffy14rates}].

In this respect, the availability of Theorem~\ref{thm:generalIP}
allows us
to improve upon the analysis provided in an earlier paper by
\citet{CW13mimeo}---and, in the random walk case, that of \citet{GLT13mimeo}---who
obtain uniform convergence rates for $\hat{m}_{n}$ on precisely such
intervals [see \citet{Duffy14rates} for further details]. We expect
that it would also play a similarly important role in the derivation
of uniform convergence rates for series regression estimators in this
setting, by ensuring the eigenvalues of the design matrix diverge
at an appropriate rate, when attempting to estimate $m_{0}$ on a
sequence of domains that contains most of the observed $\{x_{t}\}_{t=1}^{n}$.

\section{Proof of Theorem~\texorpdfstring{\protect\ref{thm:generalIP}}{3.1}}
\label{sec:IPproof}

We shall prove only part~(i) of Theorem~\ref
{thm:generalIP} here;
the relatively minor modifications required for the proof of
part~(ii)
are detailed in Section~B of the Supplement. Let $M<\infty$ be given;
it suffices to prove that (\ref{eq:mainIP}) holds in $\ell_{\infty}[-M,M]$.
To simplify the exposition, we shall require that $h_{n}\in\mathscr{H}_{n}$
always; the proof in the general case (where this occurs w.p.a.1)
requires no new ideas. The proof involves three steps:
\begin{longlist}[(iii)]
\item[(i)] show that $\mathcal{L}_{n}^{\varphi}(a)\rightsquigarrow
\mathcal{L}(a)$, using Proposition~\ref{prop:increment};
\item[(ii)] deduce $\mathcal{L}_{n}^{f}(a,h_{n})\rightsquigarrow\mathcal
{L}(a)\int_{\mathbb{R}}f$ for
$f\in\mathrm{BIL}_{\beta}$,
using a recentering, Proposition~\ref{prop:maxSn} and the Lipschitz continuity
of $f$;
\item[(iii)] extend this to all $f\in\mathscr{F}\subset\mathrm{BI}$, where
$\mathscr{F}$ satisfies
Assumption~\ref{ass:fset}, via a bracketing argument.
\end{longlist}

\begin{longlist}[(i)]
\item[(i)] Let $\varphi$ be the triangular kernel function, as defined in (\ref{eq:cmpct})
above, and set $\beta_{0}:=\overline{\beta}_{H}/2$. Recall that
$\mathcal{L}_{n}^{f}(a):=\mathcal{L}_{n}^{f}(a,1)$.
By Proposition~\ref{prop:increment} and Theorem~2.2.4 in \citet{VVW96},
\begin{eqnarray*}
&& \Bigl\Vert\sup_{\{a,a^{\prime}\in M\vert\vert a-a^{\prime}\vert\leq
\delta\} }\bigl\vert\mathcal{L}_{n}^{\varphi}
\bigl(a^{\prime}\bigr)-\mathcal {L}_{n}^{\varphi}(a)
\bigr\vert
\Bigr\Vert_{1}
\\
&& \qquad \lesssim\int_{0}^{\delta}
\log^{3/2}\bigl(M\varepsilon^{-1/\beta_{0}}\bigr)\, \mathrm{d}
\varepsilon+\delta\log^{3/2}\bigl(M\delta^{-2/\beta_{0}}\bigr)\lesssim
C_{M}\delta^{1/2},
\end{eqnarray*}
whence $\mathcal{L}_{n}^{\varphi}$ is tight in $\ell_{\infty
}[-M,M]$. Thus,
in view of Remark~\ref{rem:locconv},
\begin{equation}
\label{eq:wkccmpt} \mathcal{L}_{n}^{\varphi}(a)\rightsquigarrow
\mathcal{L}(a)
\end{equation}
in $\ell_{\infty}[-M,M]$ [see \citet{VVW96}, Example~2.2.12].

\item[(ii)]
Now let $f\in\mathrm{BIL}_{\beta}$; we may without loss of
generality take
$f$ to be bounded by unity, with a Lipschitz constant of unity. For
the subsequent argument, it will be more convenient to work with the
inverse bandwidth $b:=h^{-1}$. Define
\[
\mathscr{B}_{n}:=\bigl\{h^{-1}\vert h\in\mathscr{H}_{n}
\bigr\}=[\underline {b}_{n},\overline{b}_{n}]:=\bigl[
\overline{h}_{n}^{-1},\underline{h}_{n}^{-1}
\bigr]
\]
and let $f_{(a,b)}(x):=bf[b(x-d_{n}a)]$, for $(a,b)\in\mathbb
{R}\times
\mathbb{R}_{+}$.
Take $C_{n}:=[-n^{\gamma},n^{\gamma}]\times\mathscr{B}_{n}$, let
$\mathscr{C}
_{n}\subset C_{n}$
be a lattice of mesh $n^{-\delta}$, and let $p_{n}(a,b)$ denote
the projection of $(a,b)\in C_{n}$ onto a nearest neighbor in $\mathscr{C}_{n}$
(with some tie-breaking rule). The following is a straightforward
consequence of the Lipschitz continuity of $f$ (see Section~C of
the Supplement for the proof).
\end{longlist}

\begin{lem}
\label{lem:gridapprox}
For every $\gamma\geq1$, there exists $\delta>0$
such that
\[
\sup_{(a,b)\in C_{n}}\frac{1}{e_{n}}\sum_{t=1}^{n}
\bigl\vert f_{(a,b)}(x_{t})-f_{p_{n}(a,b)}(x_{t})
\bigr\vert=o_{p}(1).
\]
\end{lem}

By taking $\gamma\geq1$, we may ensure that $C_{n}\supset
[-M,M]\times
\mathscr{B}_{n}$,
for all $n$ sufficiently large. Thus, for $\varphi_{(a)}:=\varphi_{(a,1)}$
and $\mu_{f}:=\int_{\mathbb{R}}f$,
\begin{eqnarray}
\label{eq:centredLip} &&\sup_{(a,b)\in[-M,M]\times\mathscr{B}_{n}}\bigl\vert\mathcal {L}_{n}^{f}
\bigl(a,b^{-1}\bigr)-\mu _{f}\mathcal{L}_{n}^{\varphi}(a)
\bigr\vert
\\
&&\qquad \leq\sup_{(a,b)\in C_{n}}\frac{1}{e_{n}}\Biggl|\sum
_{t=1}^{n}\bigl[f_{(a,b)}(x_{t})-
\mu_{f}\varphi_{(a)}(x_{t})\bigr]\Biggr|
\nonumber
\\
&&\qquad\leq\sup_{(a,b)\in\mathscr{C}_{n}}\frac
{1}{e_{n}}\Biggl|\sum
_{t=1}^{n}\bigl[f_{(a,b)}(x_{t})-
\mu_{f}\varphi _{(a)}(x_{t})\bigr]\Biggr|+o_{p}(1)
\nonumber
\\
\label{eq:centredLipend} && \qquad=\sup_{g\in\mathscr{G}_{n}}\frac
{1}{e_{n}}\Biggl|\sum
_{t=1}^{n}g(x_{t})\Biggr|+o_{p}(1),
\end{eqnarray}
by Lemma~\ref{lem:gridapprox}, and we have\vspace*{1pt} defined $\mathscr
{G}_{n}:=\{
f_{(a,b)}-\mu_{f}\varphi_{(a)}\vert(a,b)\in\mathscr{C}_{n}\}$.
It is readily verified that $\Vert g\Vert_{1}=1$, $\#\mathscr
{G}_{n}=\#\mathscr{C}
_{n}\lesssim n^{1+\gamma+2\delta}$,
and using Lemma~\ref{lem:alphnorm}(ii),
\[
\sup_{g\in\mathscr{G}_{n}}\Vert g\Vert_{[\beta]} \lesssim
\underline{b}_{n}^{-\beta
}=o\bigl(d_{n}^{\beta}
\bigr), \qquad \sup_{g\in\mathscr{G}_{n}} \Vert g\Vert_{\infty
}\leq
\overline{b}_{n}\lesssim e_{n}\log^{-2}n.
\]
Thus, $\mathscr{G}_{n}$ satisfies the requirements of Proposition~\ref
{prop:maxSn},
whence (\ref{eq:centredLipend}) is $o_{p}(1)$. Hence, in view of (\ref{eq:wkccmpt}),
\begin{equation}
\label{eq:wkcBIL} \mathcal{L}_{n}^{f}(a,h_{n})
\rightsquigarrow\mathcal{L}(a)\int_{\mathbb{R}}f
\end{equation}
in $\ell_{\infty}[-M,M]$, for every $f\in\mathrm{BIL}_{\beta}$.

\begin{longlist}[(iii)]
\item[(iii)]
Finally, for $f\in\mathrm{BI}$ define the centered process
\[
\nu_{n}(f,a):=\mathcal{L}_{n}^{f}(a,h_{n})-
\mathcal{L}_{n}^{\varphi
}(a)\int_{\mathbb{R}}f.
\]
For a given $\varepsilon>0$, let $\{l_{k},u_{k}\}_{k=1}^{K}$ denote
a collection of \textit{continuous} $L^{1}$ brackets that cover
$\mathscr{F}$,
with $\Vert u_{k}-l_{k}\Vert_{1}<\varepsilon$; the existence of these
is guaranteed by Assumption~\ref{ass:fset}. We first note (see
Section~C of
the Supplement  for the proof)
\end{longlist}

\begin{lem}
\label{lem:Lipbkt}
Under Assumption~\ref{ass:fset}, the brackets $\{l_{k},u_{k}\}_{k=1}^{K}$
can be chosen so as to lie in $\mathrm{BIL}_{\beta}$.
\end{lem}

For each $f\in\mathscr{F}$, there exists a $k\in\{1,\ldots,K\}$
such that
$l_{k}\leq f\leq u_{k}$, $\int_{\mathbb{R}}(u_{k}-f)<\varepsilon$, and
\begin{eqnarray*}
\nu_{n}(f,a) & \leq & \frac{1}{e_{n}}\sum
_{t=1}^{n} \biggl[\frac
{1}{h_{n}}u_{k}
\biggl(\frac{x_{t}-d_{n}x}{h_{n}} \biggr)-\varphi (x_{t})\int
_{\mathbb{R}}f \biggr]
\\
& \leq & \nu_{n}(u_{k},a)+\mathcal{L}_{n}^{\varphi}(a)
\int_{\mathbb
{R}}(u_{k}-f).
\end{eqnarray*}
Taking suprema,
\begin{eqnarray*}
\sup_{(f,a)\in\mathscr{F}\times[-M,M]}\nu_{n}(f,a) & \leq & \max
_{1\leq k\leq
K}\sup_{a\in[-M,M]}\nu_{n}(u_{k},a)+
\varepsilon\sup_{a\in
[-M,M]}\mathcal{L}_{n}
^{\varphi}(a)
\\
& =& \varepsilon\sup_{a\in[-M,M]}\mathcal{L}_{n}^{\varphi}(a)+o_{p}(1)
\end{eqnarray*}
with the second equality following by (\ref{eq:wkcBIL}), since we may
take $u_{k}\in\mathrm{BIL}_{\beta}$ by Lemma~\ref{lem:Lipbkt}.
Applying a strictly
analogous argument to the lower bracketing functions, $l_{k}$, we
deduce that
\begin{equation}
\label{eq:bktnegl} \sup_{(f,a)\in\mathscr{F}\times[-M,M]}\bigl\vert\nu_{n}(f,a)\bigr\vert\leq
\varepsilon\sup_{a\in[-M,M]}\bigl\vert\mathcal{L}_{n}^{\varphi}(a)
\bigr\vert +o_{p}(1)=o_{p}(1),
\end{equation}
whence (\ref{eq:mainIP}) holds in $\ell_{\infty}[-M,M]$, in view of
(\ref{eq:wkccmpt}).

\section{Controlling the zero energy process}
\label{sec:ptwise}

The proofs of Propositions~\ref{prop:increment} and \ref{prop:maxSn}
rely on a telescoping martingale decomposition similar to that used
to prove maximal inequalities for mixingales [for a textbook exposition,
see, e.g.,  \citet{Davidson94}, Sections~16.2--16.3], which reduces
$\S_{n}f$ to a sum of martingale components. In order to pass from
control over each of these components to an order estimate for $\S_{n}f$
itself, we shall need the following results, the first of which is
a straightforward consequence of Theorem~2.1 in \citet{BT08AAP},
and the second of which is well known. For a martingale $M:=\{
M_{t}\}_{t=0}^{n}$
with associated filtration $\mathcal{G}:=\{\mathcal{G}_{t}\}_{t=0}^{n}$,
define
\begin{equation}
\label{eq:qvcv} [M] := \sum_{t=1}^{n}(M_{t}-M_{t-1})^{2},\qquad
\langle M\rangle := \sum_{t=1}^{n}
\mathbb{E}\bigl[(M_{t}-M_{t-1})^{2}\vert\mathcal
{G}_{t-1}\bigr].
\end{equation}
We say that $M$ is \textit{initialised at zero} if $M_{0}=0$. Let
$\Vert\cdot\Vert_{\tau_{1}}$ denote the Orlicz norm associated to
$\tau_{1}(x):=\mathrm{e}^{x}-1$.

\begin{lem}
\label{lem:mgmax} Let $\{\Theta_{n}\}$ denote a sequence of index
sets, and $\{K_{n}\}$ a real sequence such that $\#\Theta
_{n}+K_{n}\lesssim n^{C}$.
Suppose that for each $n\in\mathbb{N}$, $k\in\{1,\ldots,K_{n}\}$
and $\theta\in\Theta_{n}$, $M_{nk}(\theta)$ is a martingale, initialised
at zero, for which
\begin{equation}
\label{eq:omnk} \omega_{nk}^{2}:=\max_{\theta\in\Theta_{n}}
\bigl\{\bigl\Vert\bigl[M_{nk}(\theta )\bigr]\bigr\Vert_{\tau_{1}}\vee\bigl\Vert
\bigl\langle M_{nk}(\theta)\bigr\rangle\bigr\Vert _{\tau
_{1}}\bigr\}<
\infty.
\end{equation}
Then
\[
\max_{\theta\in\Theta_{n}}\Biggl|\sum_{k=1}^{K_{n}}M_{nk}(
\theta )\Biggr|\lesssim_{p} \Biggl(\sum_{k=1}^{K_{n}}
\omega_{nk} \Biggr)\log n.
\]
\end{lem}

\begin{lem}
\label{lem:momOrl}Let $Z$ be a random variable. Then:
\begin{longlist}[(ii)]
\item[(i)] $\Vert Z\Vert_{p}\lesssim p!^{1/p}\sigma$ for all
$p\in\mathbb{N}$, if and\vspace*{2pt} only if $\Vert Z\Vert_{\tau_{1}}\lesssim
\sigma$;
\item[(ii)] $\Vert Z\Vert_{2p}\lesssim(3p)!^{1/2p}\sigma
$ for
all $p\in\mathbb{N}$, if and only if $\Vert Z\Vert_{\tau
_{2/3}}\lesssim
\sigma$.
\end{longlist}
\end{lem}

The proofs of Lemmas~\ref{lem:mgmax} and \ref{lem:momOrl} appear
in Section~D of the Supplement.

\subsection{The martingale decomposition}
\label{sub:mgdecomp}

For a fixed $f\in\mathrm{BI}_{[\beta]}$, it follows from Lemma~\ref{lem:1stmom}(ii)
below and the reverse martingale convergence theorem [\citet{HH80},
Theorem~2.6],
that
\[
\bigl\Vert\mathbb{E}_{t}f(x_{t+k})\bigr\Vert_{\infty} \lesssim
d_{k}^{-(1+\beta
)}\rightarrow0, \qquad \mathbb{E}_{t-k}f(x_{t})
\mathop{\rightarrow }^{p}\mathbb{E}f(x_{t})\neq0
\]
for each $t\geq0$ as $k\rightarrow\infty$; here $\mathbb{E}_{t}f(x_{t+k}):=
\mathbb{E}[f(x_{t+k})\vert\mathcal{F}_{-\infty}^{t}]$,
for $\mathcal{F}_{s}^{t}:=\sigma(\{\varepsilon_{r}\}_{r=s}^{t})$. Because
$\{f(x_{t})\}$ is asymptotically unpredictable only in the ``forwards''
direction, we truncate the ``usual'' decomposition at $t=0$, writing
\[
f(x_{t}) = \sum_{k=1}^{t}\bigl[
\mathbb{E}_{t-k+1}f(x_{t})-\mathbb{E} _{t-k}f(x_{t})
\bigr]+\mathbb{E}_{0}f(x_{t}).
\]
Performing this for each $1\leq t\leq n$ gives
\begin{eqnarray*}
\sum_{t=1}^{n}f(x_{t})&=&
\mathbb{E}_{0}f(x_{1}) +\bigl[f(x_{1})-\mathbb{E}
_{0}f(x_{1})\bigr]
\\
&&{}+\mathbb{E}_{0}f(x_{2}) +\bigl[f(x_{2})-
\mathbb{E}_{1}f(x_{2})\bigr] +\bigl[\mathbb{E}_{1}f(x_{2})-
\mathbb{E}_{0}f(x_{2})\bigr]+\cdots
\\
&&{}+\mathbb{E}_{0}f(x_{n}) +\bigl[f(x_{n})-
\mathbb{E}_{n-1}f(x_{n})\bigr] +\bigl[\mathbb{E}_{n-1}f(x_{n})-
\mathbb{E}_{n-2}f(x_{n})\bigr]
\\
&&{}+\cdots+\bigl[\mathbb{E}_{1}f(x_{n})-
\mathbb{E}_{0}f(x_{n})\bigr].
\end{eqnarray*}
Defining
\begin{equation}\label{eq:xidef}
\xi_{kt}f:=\mathbb{E}_{t}f(x_{t+k})-
\mathbb{E}_{t-1}f(x_{t+k})
\end{equation}
and collecting terms appearing in the same ``column'' of the preceding
display, we thus obtain
\begin{eqnarray}
\S_{n}f & =& \sum_{t=1}^{n}f(x_{t})
=\sum_{t=1}^{n}\mathbb {E}_{0}f(x_{t})+
\sum_{k=0}^{n-1}\sum
_{t=k+1}^{n}\bigl[\mathbb{E}_{t-k}f(x_{t})-
\mathbb{E} _{t-k-1}f(x_{t})\bigr]
\nonumber
\\
& =&\sum_{t=1}^{n}\mathbb{E}_{0}f(x_{t})+
\sum_{k=0}^{n-1}\sum
_{t=1}^{n-k}\bigl[\mathbb{E}_{t}f(x_{t+k})-
\mathbb {E}_{t-1}f(x_{t+k})\bigr]
\nonumber
\\[-8pt]
\label{eq:SnDecomp}
\\[-8pt]
\nonumber
& =& \sum_{t=1}^{n}
\mathbb{E}_{0}f(x_{t})+\sum_{k=0}^{n-1}
\sum_{t=1}^{n-k}\xi_{kt}f
\\
& =& \mathcal{N}_{n}f+\sum_{k=0}^{n-1}
\mathcal{M}_{nk}f,\nonumber
\end{eqnarray}
where
\[
\mathcal{N}_{n}f :=\sum_{t=1}^{n}
\mathbb{E}_{0}f(x_{t}),\qquad \mathcal {M}_{nk}f
:= \sum_{t=1}^{n-k}\xi_{kt}f.
\]
A bound for $\Vert\mathcal{N}_{n}f\Vert_{\infty}$ is provided by
Lemma~\ref{lem:1stmom}(ii)
below. $\{\xi_{kt}f,\mathcal{F}_{-\infty}^{t}\}_{t=1}^{n-k}$, by
construction,
forms a martingale difference sequence for each $k$, and so control
over each of the martingale ``pieces'' $\mathcal{M}_{nk}f$ will follow from
control over
\[
\mathcal{U}_{nk}f :=[\mathcal{M}_{nk}f]=\sum
_{t=1}^{n-k}\xi _{kt}^{2}f,
\qquad \mathcal{V} _{nk}f :=\langle\mathcal{M}_{nk}f\rangle=
\sum_{t=1}^{n-k}\mathbb {E}_{t-1}
\xi_{kt}^{2}f,
\]
in combination with either Burkholder's inequality [\citet{HH80},
Theorem~2.10]
or Lemma~\ref{lem:mgmax} above, as appropriate.

\subsection{Proofs of Propositions~\texorpdfstring{\protect\ref{prop:increment}}{4.1} and \texorpdfstring{\protect\ref{prop:maxSn}}{4.2}}
\label{sub:propproofs}

Define
\[
\varsigma_{n}(\beta,f):=\Vert f\Vert_{\infty}+\Vert f\Vert
_{1}+\Vert f\Vert_{[\beta]}\sum_{t=1}^{n}d_{t}^{-(1+\beta)}
\]
and
\[
\sigma_{nk}^{2}(\beta,f) := \cases{ \Vert f
\Vert_{\infty}^{2}+\Vert f\Vert_{2}^{2}e_{n},
\qquad\mbox {if }k\in\{ 0,\ldots,k_{0}\},\vspace*{3pt}
\cr
e_{n} \bigl[k^{-1}d_{k}^{-(1+2\beta)}\Vert f
\Vert_{[\beta
]}^{2}+\mathrm{e} ^{-\gamma_{1}k}\Vert f
\Vert_{1}^{2} \bigr],\vspace*{3pt}\cr
\hspace*{79pt}\qquad\mbox{if }k\in
\{k_{0} +1,\ldots,n-1\}.}
\]
The following provides the requisite control over the components of
(\ref{eq:SnDecomp}).

\begin{lem}
\label{lem:fundamental}
For any $\beta\in[0,1]$,
\begin{equation}
\label{eq:Nnbnd} \Vert\mathcal{N}_{n}f\Vert_{\infty}\lesssim
\varsigma_{n}(\beta ,f),
\end{equation}
and for all $0\leq k\leq n-1$,
\begin{equation}
\label{eq:sqbnds} \Vert\mathcal{U}_{nk}f\Vert_{\tau_{1}}\vee\Vert
\mathcal {V}_{nk}f\Vert_{\tau_{1}}\lesssim \sigma_{nk}^{2}(
\beta,f).
\end{equation}
\end{lem}

The proof of (\ref{eq:sqbnds}), in turn, relies upon the following.

\begin{lem}
\label{lem:mgbnd}
For every $k\in\{0,\ldots,n-1\}$, $t\in\{1,\ldots
,n-k\}$
and $\beta\in(0,1]$
\[
\bigl\Vert\xi_{kt}^{2}f\bigr\Vert_{\infty}+\sum
_{s=1}^{n-k-t}\bigl\Vert\mathbb {E} _{t}
\xi_{k,t+s}^{2}f\bigr\Vert_{\infty}\lesssim\sigma
_{nk}^{2}(\beta,f).
\]
\end{lem}

The proofs of these results are deferred to Sections~\ref{sec:prelim}
and \ref{sec:bndaux}. We shall also need the following, for which
we recall the definition of $\delta_{n}(\beta,\mathscr{G})$ given in
(\ref{eq:delnbetaf})
above.

\begin{lem}
\label{lem:dnsum} If $\beta\in(0,\overline{\beta}_{H})$ and
$\mathscr{G}\subset
\mathrm{BI}_{[\beta]}$,
then there exists a $C_{\beta}<\infty$ such that
\[
\sup_{f\in\mathscr{G}}\varsigma_{n}(\beta,f)+\sum
_{k=0}^{n-1}\sup_{f\in\mathscr{G}
}
\sigma_{nk}(\beta,f)\leq C_{\beta}\delta_{n}(\beta,
\mathscr{G}).
\]
\end{lem}

The proof appears in Section~D of the Supplement. We now turn\vspace*{-2pt} to:

\begin{pf*}{Proof of Proposition~\ref{prop:increment}}
Let $g\in\mathrm{BI}_{[\beta]}$. Burkholder's inequality, and
Lemmas~\ref{lem:momOrl}(i)
and \ref{lem:fundamental} and give
\[
\Vert\mathcal{M}_{nk}g\Vert_{2p} \leq b_{2p}^{1/2p}
\Vert\mathcal {U} _{nk}g\Vert_{p}^{1/2}
\lesssim(b_{2p}\cdot p!)^{1/2p}\sigma _{nk}(\beta ,g)
\lesssim(3p)!^{1/2p}\sigma_{nk}(\beta,g)
\]
for every $p\in\mathbb{N}$, where $b_{2p}$ depends on $p$ in the
manner prescribed by Burkholder's inequality. Hence, $\Vert\mathcal
{M} _{nk}g\Vert_{\tau_{2/3}}\lesssim\sigma_{nk}(\beta,g)$
by Lemma~\ref{lem:momOrl}(ii). Then by (\ref
{eq:SnDecomp}), Lemma~\ref{lem:fundamental}
and Lemma~\ref{lem:dnsum} (taking\vspace*{-3pt} $\mathscr{G}=\{g\}$)
\begin{equation}
\label{eq:Sntau} \Vert\S_{n}g\Vert_{\tau_{2/3}} \leq\Vert
\mathcal{N}_{n}g\Vert _{\infty}+\sum_{k=0}^{n-1}
\Vert\mathcal{M}_{nk}g\Vert_{\tau_{2/3}}\leq C\delta _{n}(
\beta ,g)
\end{equation}
for some $C<\infty$ depending on $\beta$.

For $a_{1},a_{2}\in\mathbb{R}$, set $\Delta:=\vert a_{1}-a_{2}\vert$
and\vspace*{-2pt} define
\[
\varphi_{[a_{1},a_{2}]}(x):=\varphi(x-d_{n}a_{1})-
\varphi(x-d_{n}a_{2}).
\]
Let $\beta\in(0,\overline{\beta}_{H})$. Since $\varphi$ is bounded and\vspace*{-2pt}
Lipschitz,
\[
\Vert\varphi_{[a_{1},a_{2}]}\Vert_{\infty}\leq(d_{n}\Delta)
\wedge 1\leq d_{n}^{\beta}\Delta^{\beta};
\]
and further, since $\varphi$ is bounded and compactly supported,
\[
\Vert\varphi_{[a_{1},a_{2}]}\Vert_{p}\leq2\Vert\varphi
_{[a_{1},a_{2}]}\Vert_{\infty}^{\beta}\Vert\varphi^{1-\beta
}
\Vert_{p}\lesssim d_{n}^{\beta}\Delta^{\beta},
\]
for $p\in\{1,2\}$. Finally,\vspace*{-2pt} by Lemma~\ref{lem:alphnorm}(iii),
\[
\Vert\varphi_{[a_{1},a_{2}]}\Vert_{[\beta]}\lesssim d_{n}^{\beta
}
\Delta ^{\beta}.
\]
Then by (\ref{eq:Sntau}) and the definition of\vspace*{-2pt} $\delta_{n}(\beta
,\mathscr{F})$,
\begin{eqnarray*}
\bigl\Vert\mathcal{L}_{n}^{\varphi}(a_{1})-
\mathcal{L}_{n}^{\varphi
}(a_{2})\bigr\Vert_{\tau_{2/3}} &
=& e_{n}^{-1}\Vert\S_{n}\varphi_{[a_{1},a_{2}]}
\Vert_{\tau_{2/3}}
\\[-2pt]
& \leq & C\bigl(e_{n}^{-1/2}\cdot d_{n}^{\beta}
\Delta^{\beta
}+d_{n}^{-\beta}\cdot d_{n}^{\beta}
\Delta^{\beta}\bigr)
\\[-2pt]
& \lesssim & C\Delta^{\beta},
\end{eqnarray*}
for some $C$ depending on $\beta$; here we have used the fact that
since $\beta<\overline{\beta}_{H}\leq\frac{1-H}{2H}$, $\{e
_{n}^{-1/2}d_{n}^{\beta}\}$
is regularly varying with index $H(\beta-\frac{1-H}{2H})<0$.
\end{pf*}

\begin{pf*}{Proof of Proposition~\protect\ref{prop:maxSn}}
In view of Lemmas~\ref{lem:fundamental} and \ref{lem:dnsum}, we
have
\[
\max_{f\in\mathscr{F}_{n}}\vert\mathcal{N}_{n}f\vert\leq\max
_{f\in\mathscr{F}_{n}}\varsigma _{n}(\beta,f)\lesssim_{p}
\delta_{n}(\beta,\mathscr{F}_{n}),
\]
and\vspace*{-2pt} by an application of Lemma~\ref{lem:mgmax},
\[
\max_{f\in\mathscr{F}_{n}}\Biggl|\sum_{k=0}^{n-1}
\mathcal {M}_{nk}f\Biggr|\lesssim_{p}\delta _{n}(\beta,
\mathscr{F}_{n})\log n.
\]
Thus, (\ref{eq:Snbnd}) follows from (\ref{eq:SnDecomp}).

For the second part of the result, note that under the stated conditions
on $\mathscr{F}_{n}$,
\[
\Vert\mathscr{F}_{n}\Vert_{2}\leq\Vert\mathscr{F}_{n}
\Vert _{\infty}^{1/2}\Vert\mathscr{F}_{n}\Vert
_{1}^{1/2}=o\bigl[e_{n}^{1/2}
\log^{-1}n\bigr],
\]
whence
\[
e_{n}^{-1}\delta_{n}(\beta,\mathscr{F}_{n})=o_{p}
\bigl(\log ^{-1}n\bigr)+d_{n}^{-\beta}o_{p}
\bigl(d_{n}^{\beta}\bigr)=o_{p}(1)
\]
whereupon the result follows by (\ref{eq:Snbnd}).
\end{pf*}

\section{Proof of (\texorpdfstring{\protect\ref{eq:mostsupp}}{5.4})}
\label{sec:proofsupp}

Let $\mathcal{\mu}_{n}(a):=\frac{1}{n}\sum_{t=1}^{n}\mathbf{1}\{
d_{n}^{-1}x_{t}\leq a\}$
and $\mathcal{\mu}(a):=\break \int_{-\infty}^{a}\mathcal{L}(x)\,\mathrm
{d}x$. It is
shown in Section~E of the Supplement that
\begin{equation}
\mathcal{\mu}_{n}\rightsquigarrow\mathcal{\mu}\label{eq:wkclmeas}
\end{equation}
in $\ell_{\infty}(\mathbb{R})$, jointly with the convergence in
Theorem~\ref{thm:generalIP}.

Let $T(x):=\mathbf{1}\{x<\varepsilon\}$ and $\mathcal
{L}_{n}(a):=\mathcal{L}_{n}^{K}(a,h_{n})$.
We first note that
\[
\frac{1}{n}\sum_{t=1}^{n}\mathbf{1}
\bigl\{x_{t}\notin A_{n}^{\varepsilon}\bigr\} =
\frac
{1}{n}\sum_{t=1}^{n}\mathbf{1}
\bigl\{\mathcal {L}_{n}\bigl(d_{n}^{-1}x_{t}
\bigr)<\varepsilon\bigr\}=\int_{\mathbb{R}}T\bigl(
\mathcal{L}_{n}(a)\bigr)\,\mathrm{d}\mathcal{\mu}_{n}(a)
\]
by definition of $A_{n}^{\varepsilon}$ and $\mathcal{\mu}_{n}$. We
shall now suppose
that $\mathcal{L}_{n}\stackrel{\mathrm{a.s.}}{\rightarrow}\mathcal
{L}_{n}$ in $\ell_{\mathrm{ucc}}(\mathbb{R})$, and $\mathcal{\mu
}_{n}\stackrel{\mathrm{a.s.}}{\rightarrow}
\mathcal{\mu}$
in $\ell_{\infty}(\mathbb{R})$, as may be justified [in view of
Theorem~\ref{thm:generalIP}
and (\ref{eq:wkclmeas})] by Theorem~1.10.3 in \citet{VVW96}; let
$\Omega
_{0}\subset\Omega$
denote a set, having $\mathbb{P}\Omega_{0}=1$, on which this convergence
occurs. Define
\[
\overline{T}(x)= %
\cases{ 1, & \quad$\mbox{if }x\leq\varepsilon$,
\vspace*{3pt}
\cr
\varepsilon^{-1}(2\varepsilon-x), & $\quad\mbox{if }x
\in(\varepsilon ,2\varepsilon)$,\vspace*{3pt}
\cr
0, & $\quad\mbox{if }x\geq2
\varepsilon$.}
\]
Then, fixing an $\omega\in\Omega_{0}$,
\[
\int_{\mathbb{R}}T\bigl(\mathcal{L}_{n}^{\omega}(a)
\bigr)\,\mathrm {d}\mathcal{\mu}_{n}^{\omega}(a) \leq\int
_{\mathbb{R}}\overline{T}\bigl(\mathcal{L}_{n}^{\omega}(a)
\bigr)\,\mathrm {d}\mathcal{\mu}_{n}^{\omega}(a)=\int
_{\mathbb{R}
}F_{n}(a)\,\mathrm{d}\mathcal{
\mu}_{n}^{\omega}(a),
\]
where $F_{n}(a):=(\overline{T}\circ\mathcal{L}_{n}^{\omega})(a)$.
Now let
$[c,d]$ be chosen such that $\mathcal{\mu}^{\omega}(d)-\break \mathcal{\mu
}^{\omega
}(c)<\varepsilon$.
Since $\overline{T}$ is uniformly continuous, $F_{n}(a)\rightarrow F(a):=
(\overline{T}\circ\mathcal{L}^{\omega})(a)$
uniformly over $a\in[c,d]$, whence
\begin{eqnarray*}
\int_{\mathbb{R}}F_{n}(a)\,\mathrm{d}\mathcal{
\mu}_{n}^{\omega
}(a) &\leq & \int_{[c,d]^{c}}\,
\mathrm{d} \mathcal{\mu}_{n}(a)+\int_{[c,d]}F(a)\,
\mathrm{d}\mathcal{\mu }_{n}^{\omega}(a)+\sup_{a\in
[c,d]}
\bigl\vert F_{n}(a)-F(a)\bigr\vert
\\
&\rightarrow & \varepsilon+\int_{[c,d]}F(a)\,\mathrm{d}
\mathcal{\mu }^{\omega}(a)
\\
& \leq & \varepsilon+\int_{\mathbb{R}}F(a)\,\mathrm{d}\mathcal{\mu
}^{\omega}(a),
\end{eqnarray*}
where the convergence follows by the Portmanteau theorem [\citet{VVW96}, Theorem~1.3.4],
since $F$ is continuous.

Thus,
\begin{eqnarray}
\limsup_{n\rightarrow\infty}\mathbb{P} \Biggl\{ \frac{1}{n}\sum
_{t=1}^{n}\mathbf{1}\bigl\{ x_{t}
\notin A_{n}^{\varepsilon}\bigr\}\geq\delta \Biggr\} & \leq & \limsup
_{n\rightarrow
\infty}\mathbb{P} \biggl\{ \int_{\mathbb{R}}
\overline{T}\bigl(\mathcal {L}_{n}(a)\bigr)\,\mathrm{d}\mathcal{
\mu}_{n}(a)\geq \delta \biggr\}
\nonumber
\\
\label{eq:supptarget}
&\leq & \mathbb{P} \biggl\{ \varepsilon+\int_{\mathbb{R}}
\overline {T}\bigl(\mathcal{L}(a)\bigr)\,\mathrm{d}\mathcal{\mu} (a)\geq\delta
\biggr\}
\\
\nonumber
& \leq & \mathbb{P} \biggl\{ \varepsilon+\int_{\mathbb{R}}
\mathbf{1}\bigl\{ \mathcal{L}(a)\leq 2\varepsilon\bigr\}\mathcal{L}(a)\,
\mathrm{d}a\geq\delta \biggr\},
\end{eqnarray}
where noting that $\mathcal{L}$ is the density of $\mathcal{\mu}$,
the final
inequality follows from
\[
\int_{\mathbb{R}}\overline{T}\bigl(\mathcal{L}(a)\bigr)\,\mathrm{d}
\mathcal {\mu}(a) \leq\int_{\mathbb{R}}\mathbf{1} \bigl\{
\mathcal{L}(a)\leq2\varepsilon\bigr\}\,\mathrm{d}\mathcal{\mu}(a)=\int
_{\mathbb{R}}\mathbf{1}\bigl\{\mathcal{L} (a)\leq2\varepsilon\bigr\}
\mathcal{L}(a)\,\mathrm{d}a.
\]
Finally,
\[
\varepsilon+\int_{\mathbb{R}}\mathbf{1}\bigl\{\mathcal{L}(a)\leq2
\varepsilon \bigr\}\mathcal{L} (a)\,\mathrm{d}a\mathop{\rightarrow}^{\mathrm{a.s.}}
\int_{\mathbb
{R}}\mathbf{1}\bigl\{\mathcal{L}(a)=0\bigr\}
\mathcal{L}(a)\,\mathrm{d}a=0
\]
as $\varepsilon\rightarrow0$, by dominated convergence, and so $\varepsilon>0$
may be chosen such that the right side of (\ref{eq:supptarget}) is less
than $\delta$.

\section{Results preliminary to the proofs of Lemmas~\texorpdfstring{\protect\ref{lem:fundamental}}{7.3}
and \texorpdfstring{\protect\ref{lem:mgbnd}}{7.4}}
\label{sec:prelim}

Our arguments shall rely heavily on the use of the inverse Fourier
transform to analyse objects of the form $\mathbb{E}_{t}f(x_{t+k})$,
similarly to \citet{BI95Stek}, \citeauthor{Jeg04} (\citeyear{Jeg04,Jeg08CDFP}) and
\citeauthor{WP09Ecta} (\citeyear{WP09Ecta,WP11ET}).
Provided that $f\in\mathrm{BI}$ and $Y$ has an integrable characteristic
function $\psi_{Y}$, the ``usual'' inversion formula
\begin{equation}
\label{eq:fourinv} \mathbb{E}f(y_{0}+Y)=\frac{1}{2\pi}\int
_{\mathbb{R}}\hat {f}(\lambda)\mathrm{e}^{-\i
\lambda y_{0}}\mathbb{E}
\mathrm{e}^{-\i\lambda Y}\,\mathrm {d}\lambda
\end{equation}
for $y_{0}\in\mathbb{R}$, is still valid, even when $\hat{f}(\lambda
)=\int
f(x)\mathrm{e}^{\i\lambda x}\,\mathrm{d}x$
is not integrable; these conditions will always be met whenever the
inversion formula is required below. The following provides some useful
bounds for $\hat{f}$.

\begin{lem}
\label{lem:alphnorm}
For every $f\in\mathrm{BI}$ and $\beta\in(0,1]$,
\begin{longlist}[(iii)]
\item[(i)] $\vert\hat{f}(\lambda)\vert\leq(\vert
\lambda\vert^{\beta}\Vert f\Vert_{[\beta]})\wedge\Vert f\Vert_{1}$;
\item[(ii)] if $\int f=0$, then
\[
\Vert f\Vert_{[\beta]}\leq2^{1-\beta}\inf_{y\in\mathbb{R}}
\int_{\mathbb{R}}\bigl\vert f(x-y)\bigr\vert\vert x\vert^{\beta}\,
\mathrm{d}x
\]
and so $\mathrm{BI}_{[\beta]}\supseteq\{f\in\mathrm{BI}_{\beta
}\vert\int f=0\}$;
\item[(iii)] if $f(x):=g(x-a_{1})-g(x-a_{2})$ for some
$a_{1},a_{2}\in\mathbb{R}$, then
\[
\Vert f\Vert_{[\beta]}\leq2^{1-\beta}\vert a_{1}-a_{2}
\vert^{\beta
}\Vert g\Vert_{1}.
\]
\end{longlist}
\end{lem}

Let $\mathcal{F}_{s}^{t}:=\sigma(\{\varepsilon_{r}\}_{r=s}^{t})$, noting
that $\mathcal{F}_{s_{1}}^{s_{2}}\perpp\mathcal{F}_{s_{3}}^{s_{4}}$
for $s_{1}\leq
s_{2}<s_{3}\leq s_{4}$.
For $0<\break s<t$, we shall have frequent recourse to the following decomposition:
\begin{equation}\label{eq:decomp1}
\qquad x_{t} = \sum_{k=1}^{t}v_{t}=
\sum_{k=1}^{t}\sum
_{l=0}^{\infty}\phi _{l}\varepsilon_{k-l}
%\nonumber
%\\[-9pt]
%
%\\[-9pt]
%\nonumber
%&
=:  x_{s-1,t}^{\ast}+\sum
_{i=0}^{t-s}\varepsilon_{t-i}\sum
_{j=0}^{i}\phi _{j}=:x_{s-1,t}^{\ast}+x_{s,t,t}^{\prime},
\end{equation}
where $x_{s-1,t}^{\ast}\perpp x_{s,t,t}^{\prime}$ and $x_{s-1,t}^{\ast}$
is $\mathcal{F}_{-\infty}^{s-1}$-measurable.%
\footnote{$x_{s-1,t}^{\ast}$ is weighted sum of $\{\varepsilon_{t}\}
_{t=-\infty}^{s-1}$:
since these weights are not important for our purposes, we have refrained
from giving an explicit formula for these here.%
} Defining $a_{i}:=\sum_{j=0}^{i}\phi_{j}$, we may further decompose
$x_{s,t,t}^{\prime}$ as
\begin{equation}
\label{eq:decomp2} x_{s,t,t}^{\prime}=\sum
_{i=s}^{t}a_{t-i}\varepsilon_{i}=
\sum_{i=s}^{r}a_{t-i}
\varepsilon_{i}+\sum_{i=r+1}^{t}a_{t-i}
\varepsilon_{i}=: x_{s,r,t}^{\prime}+x_{r+1,t,t}^{\prime},
\end{equation}
where $x_{s,r,t}^{\prime}$ is $\mathcal{F}_{s}^{r}$-measurable, and
$x_{r+1,t,t}^{\prime}$
is $\mathcal{F}_{r+1}^{t}$-measurable. The following property of the
coefficients
$\{a_{i}\}$ is particularly important: there exist $0<\underline
{a}\leq\break \overline{a}<\infty$,
and a $k_{0}\in\mathbb{N}$ such that
\begin{equation}
\label{eq:coefcf} \underline{a}\leq\inf_{k_{0}+1\leq k}\inf
_{\lfloor k/2\rfloor\leq
l\le
k}c_{k}^{-1}\vert a_{l}
\vert\leq\sup_{k_{0}+1\leq k}\sup_{\lfloor
k/2\rfloor\leq l\le k}c_{k}^{-1}
\vert a_{l}\vert\leq\overline {a}.
\end{equation}
This is an easy consequence of Karamata's theorem. Throughout the
remainder of the paper, $k_{0}$ refers to the object of (\ref{eq:coefcf});
it is also implicitly maintained $k_{0}\geq8p_{0}$ for $p_{0}$
as in Assumption~\ref{ass:reg}(i).

Having decomposed $x_{t}$ into a sum of independent components, we
shall proceed to control such objects as the right-hand side of (\ref{eq:fourinv})
with the aid of Lemma~\ref{lem:alphnorm} and the following, which provides
bounds on integrals involving the characteristic functions of some
of those components of $x_{t}$. Recall that Assumption~\ref{ass:reg}(i)
is equivalent to the statement that
\begin{equation}
\label{eq:dofattr} \log\psi(\lambda)=-\vert\lambda\vert^{\alpha}G(\lambda)
\biggl[1+\i\beta \operatorname{sgn}(\lambda)\tan \biggl(\frac{\pi\alpha}{2} \biggr)
\biggr]
\end{equation}
for all $\lambda$ in a neighborhood of the origin, where $G$ is
even and slowly varying at zero [see \citet{IL71}, Theorem~2.6.5].
Here, as throughout the remainder of this paper, a slowly varying
(or regularly varying) function is understood to take only strictly
positive values, and have the property that $G(\lambda)=G(\vert
\lambda\vert)$
for every $\lambda\in\mathbb{R}$.

\begin{lem}
\label{lem:covbnd}
Let $p\in[0,5]$, $q\in(0,2]$ and $z_{1},z_{2}\in
\mathbb{R}_{+}$.
Then:
\begin{longlist}[(ii)]
\item[(i)] there exists a $\gamma_{1}>0$ such that, for every
$t\geq0$ and $k\geq k_{0}+1$,
\[
\int_{\mathbb{R}}\bigl(z_{1}\vert\lambda
\vert^{p}\wedge z_{2}\bigr)\bigl\vert \mathbb{E}\mathrm{e}
^{-\i\lambda x_{t+1,t+k,t+k}^{\prime}}\bigr\vert\, \mathrm{d}\lambda\lesssim z_{1}d_{k}^{-(1+p)}+z_{2}
\mathrm{e}^{-\gamma_{1}k}
\]
and if $F(u)\asymp G^{p/\alpha}(u)$ as $u\rightarrow0$,
\begin{eqnarray*}
&& \int_{\mathbb{R}}\bigl(z_{1}\vert a_{k}
\vert^{p}\vert\lambda \vert ^{p+q}F(a_{k}\lambda)
\wedge z_{2}\bigr)\bigl\vert\mathbb{E}\mathrm{e}^{-\i
\lambda x_{t+1,t+k,t+k}^{\prime}}\bigr\vert\,
\mathrm{d}\lambda
\\
&& \qquad \lesssim z_{1}k^{-p/\alpha}d_{k}^{-(1+q)}+z_{2}
\mathrm{e}^{-\gamma_{1}k};
\end{eqnarray*}

\item[(ii)] for every $t\geq1$, $k\geq k_{0}+1$ and $s\in\{
k_{0}+1,\ldots,t\}$,
\[
\int_{\mathbb{R}}\bigl\vert\mathbb{E}\mathrm{e}^{-\i\lambda
x_{t-s+1,t-1,t+k}^{\prime }}\bigr\vert\,
\mathrm{d}\lambda\lesssim\frac
{c_{s}}{c_{k+s}}d_{s}^{-1}.
\]
\end{longlist}
\end{lem}

The preceding summarizes and refines some of the calculations presented
on pages~15--21 of \citet{Jeg08CDFP}. It further implies:

\begin{lem}
\label{lem:1stmom}
Let $f\in\mathrm{BI}$. Then:
\begin{longlist}[(ii)]
\item[(i)] for every $t\geq0$ and $k\geq k_{0}+1$
\[
\mathbb{E}_{t}\bigl\vert f(x_{t+k})\bigr\vert\lesssim
d_{k}^{-1}\Vert f\Vert_{1};
\]

\item[(ii)] if in addition $f\in\mathrm{BI}_{[\beta]}$,
then for
every $t\geq0$ and $k\geq k_{0}+1$,
\[
\bigl\vert\mathbb{E}_{t}f(x_{t+k})\bigr\vert\lesssim
\mathrm{e}^{-\gamma
_{1}k}\Vert f\Vert_{1}+d_{k}^{-(1+\beta)}
\Vert f\Vert_{[\beta]}.
\]
\end{longlist}
\end{lem}

For the next result, define
\[
\vartheta(z_{1},z_{2}):=\mathbb{E} \bigl[
\mathrm{e}^{-\i
z_{1}\varepsilon
_{0}}-\mathbb{E}\mathrm{e}^{-\i z_{1}\varepsilon_{0}} \bigr] \bigl[
\mathrm{e}^{-\i
z_{2}\varepsilon_{0}}-\mathbb{E}\mathrm{e}^{-\i z_{2}\varepsilon
_{0}} \bigr].
\]

\begin{lem}
\label{lem:cfdiffprod}
Uniformly over $z_{1},z_{2}\in\mathbb{R}$,
\[
\bigl\vert\vartheta(z_{1},z_{2})\bigr\vert\lesssim\bigl[\vert
z_{1}\vert^{\alpha
}\tilde {G}(z_{1})\wedge1
\bigr]^{1/2}\bigl[\vert z_{2}\vert^{\alpha}\tilde
{G}(z_{2})\wedge1\bigr]^{1/2},
\]
where $\tilde{G}(u)\asymp G(u)$ as $u\rightarrow0$.
\end{lem}

Proofs of (\ref{eq:fourinv}), (\ref{eq:coefcf}) and the preceding lemmas
are given in Section~F of the Supplement.

\section{Proofs of Lemmas~\protect\ref{lem:fundamental} and \protect\ref{lem:mgbnd}}
\label{sec:bndaux}
\mbox{}
\begin{pf*}{Proof of Lemma~\texorpdfstring{\ref{lem:fundamental}}{7.3}} By Lemma~\texorpdfstring{\ref{lem:1stmom}}{7.4}(ii),
\begin{eqnarray*}
\vert\mathcal{N}_{n}f\vert & \leq & \sum_{t=1}^{k_{0}}
\bigl\vert\mathbb {E} _{0}f(x_{t})\bigl\vert+\sum
_{t=k_{0}+1}^{n}\bigl\vert\mathbb {E}_{0}f(x_{t})
\bigr\vert
\\
& \lesssim & \Vert f\Vert_{\infty}+\sum_{t=k_{0}+1}^{n}
\bigl[\mathrm {e}^{-\gamma
_{1}t}\Vert f\Vert_{1}+d_{t}^{-(1+\beta)}
\Vert f\Vert_{[\beta
]} \bigr],
\end{eqnarray*}
whence (\ref{eq:Nnbnd}). Regarding (\ref{eq:sqbnds}), it follows from
repeated application of the law of iterated expectations that
\begin{eqnarray}
\qquad \mathbb{E}\vert\mathcal{V}_{nk}f\vert^{p} &\leq & p!\cdot
\sum_{t_{1}=1}^{n-k}\cdots\sum
_{t_{p-1}=t_{p-2}}^{n-k}\mathbb{E} \bigl[\mathbb{E}_{t_{1}-1}
\bigl(\xi _{kt_{1}}^{2}f\bigr)\cdots\mathbb{E}_{t_{p-1}-1}
\bigl(\xi _{kt_{p-1}}^{2}f\bigr) \bigr]
\nonumber
\\[-8pt]
\label{eq:Vnkfbnd}
\\[-8pt]
\nonumber
&&{}\times \Biggl(\bigl\Vert\xi_{kt_{p-1}}^{2}f
\bigr\Vert_{\infty}+\sum_{s=1}^{n-k-t_{p-1}}\bigl\Vert
\mathbb{E}_{t_{p-1}-1}\xi _{k,t_{p-1}+s}^{2}f\bigr\Vert_{\infty}
\Biggr),
\end{eqnarray}
more details of the calculations leading to (\ref{eq:Vnkfbnd}) are given
in Section~G of the Supplement. By Lemma~\ref{lem:mgbnd},
the final term on the right is bounded by $C\sigma_{nk}^{2}(\beta,f)$.
Proceeding inductively, we thus obtain
\[
\mathbb{E}\vert\mathcal{V}_{nk}f\vert^{p}\lesssim p!\cdot
C^{p}\sigma_{nk}^{2p}(\beta,f),
\]
whence the required bound follows by Lemma~\ref{lem:momOrl}(i).
An analogous argument yields the same bound for $\mathcal{U}_{nk}f$.
\end{pf*}

\begin{pf*}{Proof of Lemma~\ref{lem:mgbnd}}
We shall obtain the required bound for $\mathbb{E}_{t}\xi_{k,t+s}^{2}f$
by providing a bound for $\mathbb{E}_{t-s}\xi_{kt}^{2}f$ (for $s\in\{
1,\ldots,t\}$)
that depends only on $k$ and $s$ (and \textit{not} $t$), separately
considering the cases where:
\begin{longlist}[(ii)]
\item[(i)] $k\in\{k_{0}+1,\ldots,n-t\}$; and
\item[(ii)] $k\in\{0,\ldots,k_{0}\}$.
\end{longlist}

\begin{longlist}[(ii)]
\item[(i)] Recall the decomposition given in (\ref{eq:decomp1}) and (\ref{eq:decomp2})
above, applied here to reduce $x_{t+k}$ to a sum of independent pieces,
\begin{eqnarray*}
x_{t+k} & =& x_{0,t+k}^{\ast}+x_{1,t-1,t+k}^{\prime
}+x_{t,t,t+k}^{\prime
}+x_{t+1,t+k,t+k}^{\prime}
\\
& =& x_{0,t+k}^{\ast}+x_{1,t-1,t+k}^{\prime}+a_{k}
\varepsilon _{t}+x_{t+1,t+k,t+k}^{\prime}
\end{eqnarray*}
with the convention that $x_{1,t-1,t+k}^{\prime}=0$ if $t=1$, so
that by Fourier inversion,
\begin{eqnarray}
\xi_{kt}f & =& \mathbb{E}_{t}f(x_{t+k})-\mathbb
{E}_{t-1}f(x_{t+k})
\nonumber
\\
\label{eq:0enxikt}& =& \frac{1}{2\pi}\int_{\mathbb{R}}\hat{f}(\lambda)
\mathrm {e}^{-\i\lambda
x_{0,t+k}^{\ast}}\mathrm{e}^{-\i\lambda x_{1,t-1,t+k}^{\prime
}} \\
\nonumber
&&\hspace*{7pt}\qquad {}\times\bigl[\mathrm{e}^{-\i\lambda a_{k}\varepsilon
_{t}}-
\mathbb{E}\mathrm{e} ^{-\i\lambda a_{k}\varepsilon_{t}} \bigr]\mathbb{E}\mathrm{e}^{-\i
\lambda
x_{t+1,t+k,t+k}^{\prime}}\,
\mathrm{d}\lambda.
\end{eqnarray}
Then
\begin{eqnarray}
\xi_{kt}^{2}f & =&\frac{1}{(2\pi)^{2}}\int\!\!\!\int
_{\mathbb
{R}^{2}}\hat{f}(\lambda _{1})\hat{f}(
\lambda_{2})\mathrm{e}^{-\i(\lambda_{1}+\lambda
_{2})x_{0,t+k}^{\ast
}}\mathrm{e}^{-\i(\lambda_{1}+\lambda_{2})x_{1,t-1,t+k}^{\prime
}}\nonumber
\\
\label{eq:0enxiktsq} && \qquad\qquad\quad{}\times \bigl[\mathrm{e}^{-\i\lambda
_{1}a_{k}\varepsilon
_{t}}-
\mathbb{E}\mathrm{e}^{-\i\lambda_{1}a_{k}\varepsilon_{t}} \bigr] \bigl[\mathrm{e}^{-\i
\lambda_{2}a_{k}\varepsilon_{t}}-
\mathbb{E}\mathrm{e}^{-\i\lambda
_{2}a_{k}\varepsilon
_{t}} \bigr]
\\
&& \qquad\qquad\quad{}\times\mathbb{E}\mathrm{e}^{-\i\lambda
_{1}x_{t+1,t+k,t+k}^{\prime}}\mathbb{E}
\mathrm{e}^{-\i\lambda
_{2}x_{t+1,t+k,t+k}^{\prime}}\,\mathrm{d}\lambda_{1}\,\mathrm {d}
\lambda_{2}.
\nonumber
\end{eqnarray}

Now suppose $s\in\{k+1,\ldots,t\}$. Taking conditional expectations
on both sides of (\ref{eq:0enxiktsq}) gives
\begin{eqnarray*}
\mathbb{E}_{t-s}\xi_{kt}^{2}f & =&
\frac{1}{(2\pi)^{2}}\int\!\!\!\int_{\mathbb{R}^{2}}\hat {f}(
\lambda_{1})\hat{f}(\lambda_{2})\mathrm{e}^{-\i(\lambda
_{1}+\lambda
_{2})x_{0,t+k}^{\ast}}
\mathrm{e}^{-\i(\lambda_{1}+\lambda
_{2})x_{1,t-s,t+k}^{\prime}}
\\
&& \qquad\qquad\quad{}\times\mathbb{E}\mathrm{e}^{-\i(\lambda
_{1}+\lambda
_{2})x_{t-s+1,t-1,t+k}^{\prime}}\cdot\vartheta(
\lambda _{1}a_{k},\lambda _{2}a_{k})
\\
&& \qquad\qquad\quad{}\times\mathbb{E}\mathrm{e}^{-\i\lambda
_{1}x_{t+1,t+k,t+k}^{\prime}}\mathbb{E}
\mathrm{e}^{-\i\lambda
_{2}x_{t+1,t+k,t+k}^{\prime}}\,\mathrm{d}\lambda_{1}\,\mathrm {d}
\lambda_{2},
\end{eqnarray*}
where we have defined
\[
\vartheta(z_{1},z_{2}):=\mathbb{E} \bigl[
\mathrm{e}^{-\i
z_{1}\varepsilon
_{0}}-\mathbb{E}\mathrm{e}^{-\i z_{1}\varepsilon_{0}} \bigr] \bigl[
\mathrm{e}^{-\i
z_{1}\varepsilon_{0}}-\mathbb{E}\mathrm{e}^{-\i z_{2}\varepsilon
_{0}} \bigr]
\]
for $z_{1},z_{2}\in\mathbb{R}$, and made the further decomposition
\[
x_{1,t-1,t+k}^{\prime}=x_{1,t-s,t+k}^{\prime
}+x_{t-s+1,t-1,t+k}^{\prime}
\]
with the convention that $x_{1,t-s,t+k}^{\prime}=0$ if $s=t$. Then,
using (\ref{eq:coefcf}), Lemma~\ref{lem:cfdiffprod} and $\vert ab\vert
\lesssim
\vert a\vert^{2}+\vert b\vert^{2}$,
we obtain
\begin{eqnarray}
\mathbb{E}_{t-s}\xi_{kt}^{2}f & \lesssim & \int\!\!\!
\!\int_{\mathbb
{R}^{2}}\bigl\vert\hat {f}(\lambda_{1})\hat{f}(
\lambda_{2})\bigr\vert
\nonumber
\\
&& \quad\hspace*{8pt}{}\times\bigl[\vert a_{k}\lambda_{1}
\vert^{\alpha}\tilde {G}(a_{k}\lambda_{1})\wedge1
\bigr]^{1/2}\bigl[\vert a_{k}\lambda_{2}\vert
^{\alpha
}\tilde{G}(a_{k}\lambda_{2})\wedge1
\bigr]^{1/2}
\nonumber
\\[-8pt]
\label{eq:fors0case}
\\[-8pt]
\nonumber
&&\quad\hspace*{8pt}{}\times\bigl\vert\mathbb{E}\mathrm{e}^{-\i(\lambda
_{1}+\lambda _{2})x_{t-s+1,t-1,t+k}^{\prime}}\bigr\vert
\nonumber
\\
&&\quad\hspace*{8pt}{}\times\bigl\vert\mathbb{E}\mathrm{e}^{-\i\lambda
_{1}x_{t+1,t+k,t+k}^{\prime}}\bigr\vert\bigl\vert\mathbb{E}
\mathrm{e}^{-\i
\lambda _{2}x_{t+1,t+k,t+k}^{\prime}}\bigr\vert\,\mathrm{d}\lambda_{1}\, \mathrm{d}
\lambda_{2}\nonumber
\\
&\lesssim & \int_{\mathbb{R}}\bigl\vert\hat{f}(
\lambda_{1})\bigr\vert ^{2}\bigl(\vert a_{k}
\vert^{\alpha}\vert\lambda_{1}\vert^{\alpha
}
\tilde{G}(a_{k}\lambda _{1})\wedge1\bigr)\bigl\vert\mathbb{E}
\mathrm{e}^{-\i\lambda
_{1}x_{t+1,t+k,t+k}^{\prime }}\bigr\vert
\nonumber
\\[-8pt]
\label{eq:0enxikt1}
\\[-8pt]
\nonumber
&& \times\int_{\mathbb{R}}\bigl\vert\mathbb{E}
\mathrm{e}^{-\i
(\lambda_{1}+\lambda _{2})x_{t-s+1,t-1,t+k}^{\prime}}\bigr\vert\,\mathrm {d}\lambda_{2}\,\mathrm{d}
\lambda _{1},
\nonumber
\end{eqnarray}
where we have appealed to symmetry (in $\lambda_{1}$ and $\lambda_{2}$)
to reduce the final bound to a single term. By a change of variables
and Lemma~\ref{lem:covbnd}(ii),
\begin{eqnarray}
&& \int_{\mathbb{R}}\bigl\vert\mathbb{E}\mathrm{e}^{-\i(\lambda
_{1}+\lambda _{2})x_{t-s+1,t-1,t+k}^{\prime}}\bigr\vert
\,\mathrm {d}\lambda_{2}
\nonumber
\\[-8pt]
\label{eq:0enxikt2}
\\[-8pt]
\nonumber
&&\qquad=\int_{\mathbb{R}}\bigl\vert\mathbb{E}
\mathrm{e}^{-\i\lambda
x_{t-s+1,t-1,t+k}^{\prime }}\bigr\vert\,\mathrm{d}\lambda\lesssim\frac
{c_{s}}{c_{k+s}}d_{s}^{-1},
\end{eqnarray}
while Lemma~\ref{lem:alphnorm}(i) and then
Lemma~\ref{lem:covbnd}(i)
give
\begin{eqnarray}
&& \int_{\mathbb{R}}\bigl\vert\hat{f}(\lambda)\bigr\vert^{2}\bigl(
\vert a_{k}\vert^{\alpha
}\vert\lambda\vert^{\alpha}
\tilde{G}(a_{k}\lambda)\wedge1\bigr)\bigl\vert \mathbb{E}
\mathrm{e}^{-\i\lambda x_{t+1,t+k,t+k}^{\prime}}\bigr\vert\, \mathrm{d}\lambda
\nonumber
\\
\label{eq:0enxikt3} && \qquad\leq\int_{\mathbb{R}}\bigl[\bigl(\vert
a_{k}\vert^{\alpha}\vert \lambda \vert^{\alpha+2\beta}
\tilde{G}(a_{k}\lambda)\Vert f\Vert _{[\beta]}^{2}
\bigr)\wedge \Vert f\Vert_{1}^{2}\bigr]\bigl\vert\mathbb{E}
\mathrm{e}^{-\i\lambda
x_{t+1,t+k,t+k}^{\prime }}\bigr\vert\,\mathrm{d}\lambda
\\
\nonumber
&& \qquad\lesssim k^{-1}d_{k}^{-(1+2\beta)}\Vert
f\Vert_{[\beta
]}^{2}+\mathrm{e} ^{-\gamma_{1}k}\Vert f
\Vert_{1}^{2}.
\end{eqnarray}
Together, (\ref{eq:0enxikt1})--(\ref{eq:0enxikt3})
yield
\begin{equation}
\label{eq:0enmkbnd1} \mathbb{E}_{t-s}\xi_{kt}^{2}f
\lesssim\frac
{c_{s}}{c_{k+s}}d_{s}^{-1}\bigl(k^{-1}d_{k}^{-(1+2\beta)}
\Vert f\Vert _{[\beta ]}^{2}+\mathrm{e}^{-\gamma_{1}k}\Vert f\Vert
_{1}^{2}\bigr).
\end{equation}

When $s\in\{1,\ldots,k\}$, (\ref{eq:fors0case}) continues to hold,
whence
\begin{eqnarray}
\mathbb{E}_{t-s}\xi_{kt}^{2}f & \lesssim &
\biggl(\int_{\mathbb
{R}}\bigl\vert\hat {f}(\lambda)\bigr\vert\bigl(\vert\lambda
\vert^{\alpha
/2}\tilde{G}^{1/2}(a_{k}\lambda )\wedge1
\bigr)\bigl\vert\mathbb{E}\mathrm{e}^{-\i\lambda
x_{t+1,t+k,t+k}^{\prime}}\bigr\vert\,\mathrm{d} \lambda
\biggr)^{2}
\nonumber
\\[-2pt]
& \lesssim & \biggl(\int_{\mathbb{R}}\bigl[\bigl(\vert
a_{k}\vert^{\alpha
/2}\vert\lambda\vert^{(\alpha/2+\beta)}\tilde
{G}^{1/2}(a_{k}\lambda)\Vert f\Vert_{[\beta]}\bigr)
\wedge\Vert f\Vert _{1}\bigr]\nonumber \\[-3pt]
&&\label{eq:0enmkbnd1b}\hspace*{77pt}\qquad\qquad{}\times\bigl\vert\mathbb{E}\mathrm{e}^{-\i\lambda
x_{t+1,t+k,t+k}^{\prime}}
\bigr\vert\,\mathrm{d}\lambda \biggr)^{2}
\\[-3pt]
& \lesssim & \bigl(k^{-1/2}d_{k}^{-(1+\beta)}
\Vert f\Vert_{[\beta
]}+\mathrm{e}^{-\gamma
_{1}k}\Vert f\Vert_{1}
\bigr)^{2}
\nonumber
\\[-2pt]
& \lesssim& d_{s}^{-1}\bigl(k^{-1}d_{k}^{-(1+2\beta)}
\Vert f\Vert_{[\beta
]}^{2}+\mathrm{e} ^{-\gamma_{1}k}\Vert f
\Vert_{1}^{2}\bigr)
\nonumber
\end{eqnarray}
by Lemmas~\ref{lem:alphnorm}(i) and \ref
{lem:covbnd}(i);
in obtaining the final result, we have used the fact that $s\leq k$
to replace a $d_{k}^{-1}$ by $d_{s}^{-1}$. Since $\{c_{k}\}$ is
regularly varying and $k\geq k_{0}+1$, it follows from Potter's
inequality [\citet{Bingham87}, Theorem~1.5.6(iii)],\vspace*{-2pt} that
\[
\sum_{s=1}^{k}d_{s}^{-1}+
\sum_{s=k+1}^{n}\frac
{c_{s}}{c_{k+s}}d_{s}^{-1}
\lesssim\sum_{s=1}^{n}d_{s}^{-1}
\lesssim nd_{n}^{-1}=e_{n},
\]
with the final bound following by Karamata's theorem. As noted above,
since the bounds (\ref{eq:0enmkbnd1}) and (\ref{eq:0enmkbnd1b}) do not
depend on $t$, they apply also to $\mathbb{E}_{t}\xi_{k,t+s}^{2}f$.
Hence, in view of the\vspace*{-2pt} preceding,
\begin{eqnarray*}
\sum_{s=1}^{n-k-t}\mathbb{E}_{t}
\xi_{k,t+s}^{2}f & \lesssim & \bigl(k^{-1}d_{k}^{-(1+2\beta)}
\Vert f\Vert_{[\beta]}^{2}+\mathrm {e}^{-\gamma
_{1}k}\Vert f
\Vert_{1}^{2}\bigr) \Biggl[\sum
_{s=1}^{k}d_{s}^{-1}+\!\sum
_{s=k+1}^{n-k-t}\frac{c_{s}}{c_{k+s}}d_{s}^{-1}
\!\Biggr]
\\[-2pt]
& \lesssim & e_{n}\bigl(k^{-1}d_{k}^{-(1+2\beta)}
\Vert f\Vert_{[\beta
]}^{2}+\mathrm{e} ^{-\gamma_{1}k}\Vert f
\Vert_{1}^{2}\bigr).
\end{eqnarray*}

Turning now to $\Vert\xi_{kt}^{2}f\Vert_{\infty}$, note that (\ref
{eq:0enxikt})
still holds, with the convention that $x_{1,t-1,t+k}=0$ if $t=1$.
Thus, again by Lemmas~\ref{lem:alphnorm}(i) and\vspace*{-2pt}
\ref{lem:covbnd}(i),
\begin{eqnarray*}
\bigl\Vert\xi_{kt}^{2}f\bigr\Vert_{\infty} & \lesssim &
\biggl(\int_{\mathbb
{R}}\bigl\vert\hat{f}(\lambda)\bigr\vert\bigl\vert\mathbb{E}
\mathrm{e}^{-\i
\lambda x_{t+1,t+k,t+k}^{\prime }}\bigr\vert\,\mathrm{d}\lambda \biggr)^{2}
\\[-2pt]
& \lesssim & \biggl(\int_{\mathbb{R}}\bigl[\vert\lambda
\vert^{\beta
}\Vert f\Vert_{[\beta ]}\wedge\Vert f\Vert_{1}
\bigr]\bigl\vert\mathbb{E}\mathrm{e}^{-\i\lambda x_{t+1,t+k,t+k}^{\prime}}\bigr\vert\,\mathrm {d}\lambda
\biggr)^{2}
\\[-2pt]
& \lesssim & \bigl(\Vert f\Vert_{[\beta]}d_{k}^{-(1+\beta)}+
\Vert f\Vert_{1}\mathrm{e} ^{-\gamma_{1}k} \bigr)^{2}
\\[-2pt]
& \lesssim & d_{k}^{-2(1+\beta)}\Vert f\Vert_{[\beta]}^{2}+
\mathrm {e}^{-\gamma
_{1}k}\Vert f\Vert_{1}^{2}
\\[-2pt]
& \lesssim & e_{n}\bigl(k^{-1}d_{k}^{-(1+2\beta)}
\Vert f\Vert_{[\beta
]}^{2}+\mathrm{e} ^{-\gamma_{1}k}\Vert f
\Vert_{1}^{2}\bigr);
\end{eqnarray*}
where the final bound follows because $k\leq n$, and so
$d_{k}^{-1}\lesssim k^{-1}nd_{n}^{-1}=k^{-1}e_{n}$.

\item[(ii)]
When $s\in\{1,\ldots,k_{0}\}$, the crude bound $\mathbb{E}_{t-s}\xi
_{kt}^{2}f\lesssim\Vert f\Vert_{\infty}^{2}$
suffices, since $k_{0}$ is fixed and finite. On the other hand,
if $s\in\{k_{0}+1,\ldots,t\}$, we have by Jensen's inequality
and Lemma~\ref{lem:1stmom}(i) that
\[
\mathbb{E}_{t-s}\xi_{kt}^{2}f \leq
\mathbb{E}_{t-s} \bigl(\mathbb{E} _{t}f(x_{t+k})-
\mathbb{E}_{t-1}f(x_{t+k}) \bigr)^{2}\lesssim
\mathbb{E} _{t-s}f^{2}(x_{t+k})\lesssim
d_{s}^{-1}\Vert f\Vert_{2}^{2}.
\]
Then, by Karamata's theorem,
\begin{eqnarray*}
\sum_{s=1}^{n-k-t}\mathbb{E}_{t}
\xi_{k,t+s}^{2}f & \leq & \sum_{s=1}^{k_{0}
}
\mathbb{E}_{t}\xi_{k,t+s}^{2}f+\sum
_{s=k_{0}+1}^{n-k-t}\mathbb {E}_{t}\xi
_{k,t+s}^{2}f
\\
& \lesssim & \Vert f\Vert_{\infty}^{2}+\Vert f
\Vert_{2}^{2}\sum_{s=k_{0}
+1}^{n-k-t}d_{s}^{-1}
\\
& \lesssim & \Vert f\Vert_{\infty}^{2}+\Vert f
\Vert_{2}^{2}e_{n}.
\end{eqnarray*}
Regarding $\Vert\xi_{kt}^{2}f\Vert_{\infty}$, the bound $\Vert\xi
_{kt}^{2}f\Vert_{\infty}\lesssim\Vert f\Vert_{\infty}^{2}$
obtains trivially.\quad\qed
\end{longlist}
\noqed\end{pf*}

\section*{Acknowledgments}
The author thanks Xiaohong Chen, Bent Nielsen and Peter Phillips for
helpful comments and advice. He also thanks the Editor and two anonymous
referees for comments that greatly improved an earlier version of
this paper. The manuscript was prepared in L$_{\mbox{Y}}$X 2.1.2.

\bibliographystyle{imsart-nameyear}
%\bibliography{time-series}

%
% imsref loaded by daiva.urboniene, 2015-02-06 15:43:25
%

%\begin{appendix}
%\section{}
%\end{appendix}

%\section*{Acknowledgments}

%\begin{supplement}[id=suppA]
%\slink[doi]{10.1214/00-AAPXXXXSUPP} %[doi,text={...}] - jei reikia
%suskaldyti doi
%\sfilename{aapXXXX\_supp.pdf}
%\end{supplement}

%\begin{thebibliography}{99}
%\bibitem[\protect\citeauthoryear{}{}]{r1}
%\bibitem{r1}
%\end{thebibliography}

\printaddresses
\end{document}